\documentclass[a4paper,12pt]{article}

\usepackage[utf8]{inputenc}  
\usepackage[T1]{fontenc}     
\usepackage[english]{babel}  

\usepackage{geometry}    
\usepackage{graphicx} 
\usepackage{amsmath}
\usepackage{amssymb}
\usepackage{amsthm}
\usepackage{mathtools}
\usepackage[dvipsnames]{xcolor}
\usepackage{enumitem}
\usepackage{pifont}
\usepackage{hyperref}    
\usepackage{xcolor}

\newtheorem{theorem}{Theorem}[section]
\newtheorem*{theorem*}{Main Theorem}
\newtheorem*{theoremA}{Theorem A}
\newtheorem*{theoremB}{Theorem B}
\newtheorem{prop}[theorem]{Proposition}
\newtheorem{lemma}[theorem]{Lemma}
\newtheorem{coroll}[theorem]{Corollary}
\theoremstyle{definition}

\newtheorem{remark}[theorem]{Remark}

\usepackage{tikz-cd}
\usepackage{chronology}
\usepackage{epic}
\usepackage{pict2e}
\usepackage{ulem}
\newcommand{\redsout}{\bgroup\markoverwith{\textcolor{red}{\rule[0.5ex]{2pt}{2pt}}}\ULon}

\begin{document}

\begin{center}
\bf  \huge{Finite groups with a large normalized sum of element orders} 
\end{center}

\smallskip

\begin{center}
{\sc Luigi Iorio -- Marco Trombetti}  \\ \ \\
\end{center}

\begin{abstract}
For a finite group $G$, let $\psi(G)$ be the sum of the orders of its elements, and define the corresponding \textit{normalized sum} as \hbox{$\psi'(G) \coloneqq \psi(G)/\psi(\mathcal{C}_{|G|})$,} where $\mathcal{C}_{|G|}$ is the cyclic group of the same order as $G$. Inspired by analogous criteria for the classes of soluble, supersoluble, and nilpotent groups, our main result establishes that if~\hbox{$\psi'(G)>\psi'(D_8) = \frac{19}{43}$,} then $G$ belongs to the well-understood class of groups with a modular subgroup lattice, whose structure theory allows us to readily identify all groups satisfying this bound. Moreover, the equality case is fully settled. Finally, our arguments lead to a complete description of all groups satisfying \hbox{$\psi'(G)> \psi'(A_4) = \frac{31}{77}$,} thereby fully determining the groups covered by the supersolubility criterion of Baniasad Azad and Khosravi [{\it Canad. Math. Bull.} 65 (2022), 30--38], and thus providing a more complete answer to a corresponding conjecture of~T\v{a}rn\v{a}uceanu.
\end{abstract}

\medskip\medskip


\bigskip
\noindent {\bf 2020 Mathematics Subject Classification}\\
{\it Primary}: 20D60, 20E34
{\it Secondary}: 20D30, 20F16

\section{Introduction}

Let $G$ be a finite group of order $n$. The {\it $\psi$-value} of $G$ is defined as the sum of the orders of its elements, that is,
$$\displaystyle \psi(G) = \sum_{x\in G} o(x) \; ,$$
where $o(x)$ denotes the order of the element $x \in G$. The function $\psi$ was introduced in~\cite{FIRST} with the aim of investigating what information could be recovered from its value when the order of the group is fixed. The main result of~\cite{FIRST} establishes that~$\psi(G)$ is always less than or equal to the $\psi$-value of the cyclic group $\mathcal C_n$ of order $n$ (i.e., the cyclic group of the same order as $G$) with equality holding only when $G$ is cyclic. As a consequence, a finite cyclic group is uniquely determined by its order and the sum of the orders of its elements. Despite this result, a group cannot always be identified solely from its order and its $\psi$-value. For example, the two non-isomorphic semidirect products of a cyclic group of order $4$ by a cyclic group of order $4$ share the same $\psi$-value (namely, $55$). Moreover, arbitrarily large examples can be obtained by adjoining a cyclic direct factor of odd order, since the $\psi$-value of a direct product of coprime groups equals the product of the~$\psi$-values of the factors.  

To present the subsequent developments on the topic, it is convenient to define the function $\psi'$: the {\it $\psi'$-value} of $G$, which we shall also call the {\it normalized $\psi$-value} of $G$, is defined by
 $$\displaystyle \psi'(G) = \frac{\psi(G)}{\psi(\mathcal{C}_n)} \; .$$
 With this notation, the aforementioned result from $\cite{FIRST}$ can be rephrased in this way: we always have $\psi'(G) \leq 1$, and equality holds if and only if $G$ is cyclic. This fact was sharpened in \cite{NonCyclicCriterionI} by showing that if $G$ is non-cyclic, then $\psi'(G)\leq \frac{7}{11}$. Here, the constant~$\frac{7}{11}$ is not random at all: indeed, it is the $\psi'$-value of the smallest non-cyclic group, more precisely \mbox{$\psi'(\mathcal{C}_2\times\mathcal{C}_2) = \psi(\mathcal{C}_2\times\mathcal{C}_2)/\psi(\mathcal{C}_4)$.}
 This led to a shift in perspective. Indeed, the result in \cite{NonCyclicCriterionI} can be reformulated as a ``cyclicity criterion'': if the \hbox{$\psi'$-value} of $G$ is strictly greater than the \hbox{$\psi'$-value} of the smallest \hbox{non-cyclic} group, then $G$ must be cyclic.

In the wake of this, it was first conjectured in \cite{GeneralizedBoundPsiCyclic}, and then proved in \cite{SolvableCriterion}, that if~\hbox{$\psi'(G)>\frac{211}{1617}$,} then $G$ is soluble --- here, $\frac{211}{1617}$ is the $\psi'$-value of the alternating group $A_5$ of degree $5$, the smallest non-soluble group (moreover, by \cite{uguaglianza}, for non-soluble groups such a $\psi'$-value is only achieved by a direct product of $A_5$ with a coprime cyclic group). Some years later, analogous criteria have been proved for nilpotency~\cite{NilpotentCriterion} and supersolubility~\cite{SSCriterion}; here we just observe that the~\hbox{$\psi'$-values} of the smallest non-nilpotent group (\hbox{i.e.}, the symmetric group~$S_3$ of degree~$3$) and the smallest non-supersoluble group (\hbox{i.e.}, the alternating group~$A_4$ of degree~$4$) are respectively~$\frac{13}{21}$ and~$\frac{31}{77}$. All these criteria share the same form: a suitable class of groups $\mathfrak{X}$ is considered, and if $G$ has a $\psi'$-value strictly greater than the $\psi'$-value of the smallest \hbox{non-$\mathfrak{X}$ group,} then $G$ is an $\mathfrak{X}$-group. However, this does not work for every ``reasonable'' class of groups: indeed, $S_3$ is the smallest \hbox{non-abelian} group \mbox{but~$\psi'(Q_8) = \frac{27}{43} > \frac{13}{21} = \psi'(S_3)$.}
Actually, the author of \cite{NilpotentCriterion} showed that  groups whose $\psi'$-value exceeds $\frac{13}{21}$ can only be of three distinct types: cyclic groups, direct products of a Klein four group $\mathcal C_2\times\mathcal C_2$ with a cyclic group of odd order (whose $\psi'$-values are $\frac{7}{11}$), and direct products of a quaternion group~$Q_8$ of order $8$ with a cyclic group of odd order (whose $\psi'$-values are $\frac{27}{43}$); moreover, they proved that the only groups whose $\psi'$-value is~$\frac{13}{21}$ are direct products of $S_3$ with a cyclic group of coprime order. This further shifts the focus from the search for criteria concerning certain properties to the classification of groups with a large normalized sum of element orders.

The state of the art can be effectively summarized by the following diagram.

\vspace{-0.3cm}

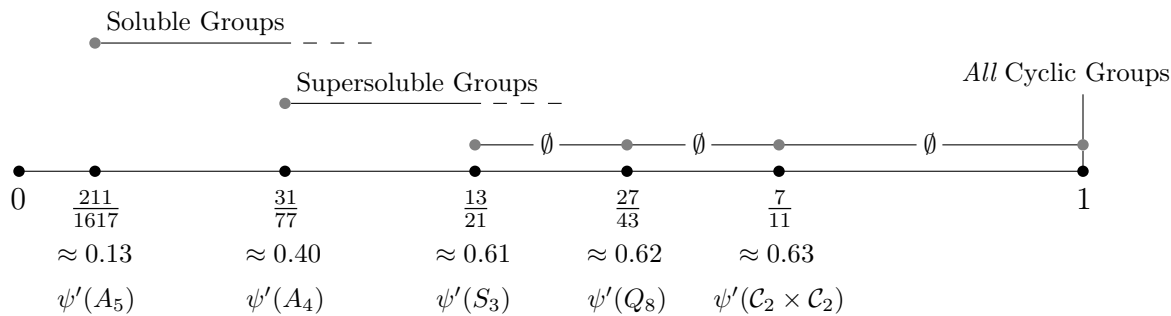
\begin{center}
\setlength{\unitlength}{1cm}
\begin{picture}(14,4)

\put(0,1){\line(1,0){14}}


\put(-0.11,0.5){0}
\put(0,1){\circle*{0.15}}
\put(0.67,0.4){$\frac{211}{1617}$}
\put(0.5,-0.2){\footnotesize $\approx 0.13$}
\put(0.5,-0.8){\footnotesize $\psi'(A_5)$}
\put(1,1){\circle*{0.15}}
\put(3.32,0.4){$\frac{31}{77}$}
\put(2.98,-0.2){\footnotesize $\approx 0.40$}
\put(3.5,1){\circle*{0.15}}
\put(2.99,-0.8){\footnotesize $\psi'(A_4)$}
\put(5.82,0.4){$\frac{13}{21}$}
\put(5.48,-0.2){\footnotesize $\approx 0.61$}
\put(5.49,-0.8){\footnotesize $\psi'(S_3)$}
\put(6,1){\circle*{0.15}}
\put(7.82,0.4){$\frac{27}{43}$}
\put(7.47,-0.2){\footnotesize $\approx 0.62$}
\put(8,1){\circle*{0.15}}
\put(7.48,-0.8){\footnotesize $\psi'(Q_8)$}
\put(9.82,0.4){$\frac{7}{11}$}
\put(9.47,-0.2){\footnotesize $\approx 0.63$}
\put(10,1){\circle*{0.15}}
\put(9.14,-0.8){\footnotesize $\psi'(\mathcal C_2\times\mathcal C_2)$}
\put(13.91,0.5){1}
\put(14,1){\circle*{0.15}}

\put(3.58,1.9){\line(1,0){2.5}}    
\put(6.24,1.9){\dashline{0.2}(0,0)(0.9,0)}
\put(3.5,1.9){\color{gray}\circle*{0.15}}
\put(3.65,2.07){\footnotesize Supersoluble Groups}

\put(1.08,2.7){\line(1,0){2.5}}   
\put(3.74,2.7){\dashline{0.2}(0,0)(0.9,0)}
\put(1,2.7){\color{gray}\circle*{0.15}}
\put(1.15,2.87){\footnotesize Soluble Groups}

\put(14,1){\line(0,1){0.27}}
\put(14,1.42){\line(0,1){0.6}}
\put(12.4,2.2){\footnotesize {\it All} Cyclic Groups}

\put(10.05,1.35){\line(1,0){1.76}}
\put(13.95,1.35){\line(-1,0){1.8}}
\put(14,1.35){\color{gray}\circle*{0.15}}
\put(11.9,1.24){\footnotesize $\emptyset$}

\put(8.05,1.35){\line(1,0){0.72}}
\put(9.95,1.35){\line(-1,0){0.84}}
\put(10,1.35){\color{gray}\circle*{0.15}}
\put(8.86,1.24){\footnotesize$\emptyset$}

\put(6.05,1.35){\line(1,0){0.72}}
\put(7.95,1.35){\line(-1,0){0.84}}
\put(8,1.35){\color{gray}\circle*{0.15}}
\put(6,1.35){\color{gray}\circle*{0.15}}
\put(6.86,1.24){\footnotesize$\emptyset$}

\end{picture}

\end{center}

\vspace{1.3cm}

As we can see, the right part of the line is rather empty, which explains why a characterization of those groups was possible in the first place. Moreover, one can also observe the following two facts: firstly, these groups are more than supersoluble, they are \hbox{\it $M$-groups} (that is, their subgroup lattice is modular); secondly, the $\psi'$-value of the dihedral group~$D_8$ of order $8$ (the smallest group that is not an $M$-group), namely $\frac{19}{43}$, lies between the $\psi'$-values of~$A_4$ and~$S_3$. Thus, by virtue of the already known criteria, and considering that finite $M$-groups are supersoluble, a natural question arises: is it possible to obtain an analogous criterion for $M$-groups in terms of the function $\psi'$? If so, these $M$-groups may turn out to be easy to describe (given that the structure of finite $M$-groups is known with remarkable precision \cite{Schmidt}), thereby leading to a better understanding of the right half of the above diagram. Our first main result completely solves this problem. Here, and throughout this paper, we use a natural and intuitive notation to denote certain semidirect products of an abelian group by a cyclic group, where the action is given by the inversion automorphism (see §2).

\medskip

\begin{theoremA}\label{main1}
Let $G$ be a finite group with
    $$\psi'(G) > \psi'(D_8) = \frac{19}{43} \ .$$ 
    Then $G$ is an $M$-group. More precisely, up to isomorphism, $G$ belongs to one of the following pairwise disjoint classes of groups:
     \begin{itemize}
        \item cyclic groups;
        \item $\left(\mathcal{C}_{2^{k-1}}\!\times \mathcal{C}_2\right) \times \mathcal{C}_m$, where $k\geq 2$ and $(m,2)=1$;
        \item $M(2^k) \times \mathcal{C}_m$, where $k\geq 4$ and $(m,2) =1$;
        \item $Q_8 \times \mathcal{C}_m$, where $(m,2)=1$;
        \item $\left(\mathcal{C}_3 \rtimes_{\iota} \mathcal{C}_{2^k}\right) \times \mathcal{C}_m$, where $k\geq 1$ and $(m,6)=1$;
        \item $\left(\mathcal{C}_5 \rtimes_{\iota} \mathcal{C}_2\right) \times \mathcal{C}_m$, where $(m,10)=1$;
        \item $\left(\mathcal{C}_5 \rtimes_{\iota} \mathcal{C}_4\right) \times \mathcal{C}_m$, where $(m,10)=1$.
    \end{itemize}
\end{theoremA}

\medskip

Our second main result deals with the extremal case in which the $\psi'$-value is required to be precisely that of $D_8$.

\medskip

\begin{theoremB}\label{thB}
Let $G$ be a finite group with
    $$\psi'(G) = \psi'(D_8) = \frac{19}{43} \ .$$
    Then, up to isomorphism, $G$ belongs to one of the following disjoint classes of groups:
    \begin{itemize}
        \item $D_8\times \mathcal{C}_m$, where $(m,2)=1$ \textnormal(in this case $G$ is not an $M$-group\textnormal);
        \item $\left(\mathcal{C}_7 \rtimes_\iota \mathcal{C}_2 \right) \times \mathcal{C}_m$, where $(m,14)=1$ \textnormal(in this case $G$ is an $M$-group\textnormal).
    \end{itemize}
\end{theoremB}

\medskip

Theorems A and B show that the interval $\big(\frac{31}{77}, \frac{13}{21}\big]$ in the above diagram can be refined as follows.
\pagebreak[1]

\vspace*{-1.7cm}

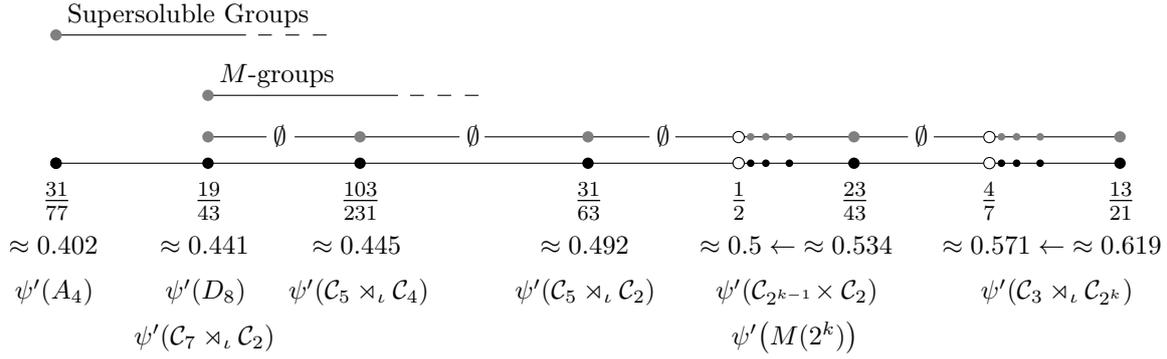
\begin{center}

\setlength{\unitlength}{1cm}
\begin{picture}(14,4)

\put(0,1){\line(1,0){14}}


\put(-0.18,0.4){$\frac{31}{77}$}
\put(-0.62,-0.2){\footnotesize $\approx 0.402$}
\put(0,1){\circle*{0.15}}
\put(-0.54,-0.8){\footnotesize $\psi'(A_4)$}
\put(0,1){\circle*{0.15}}


\put(1.82,0.4){$\frac{19}{43}$}
\put(1.37,-0.2){\footnotesize $\approx 0.441$}
\put(2,1){\circle*{0.15}}
\put(1.44,-0.8){\footnotesize $\psi'(D_8)$}
\put(1.03,-1.4){\footnotesize $\psi'(\mathcal{C}_7\rtimes_\iota\mathcal{C}_2)$}

\put(2.07,1.35){\line(1,0){0.7}}
\put(2,1.35){\color{gray}\circle*{0.15}}
\put(3.95,1.35){\line(-1,0){0.81}}
\put(2.86,1.24){\footnotesize$\emptyset$}

\put(3.74,0.4){$\frac{103}{231}$}
\put(3.37,-0.2){\footnotesize $\approx 0.445$}
\put(3.06,-0.8){\footnotesize $\psi'(\mathcal C_5\rtimes_\iota\mathcal C_4)$}
\put(4,1){\circle*{0.15}}

\put(4.05,1.35){\line(1,0){1.26}}
\put(4,1.35){\color{gray}\circle*{0.15}}
\put(6.95,1.35){\line(-1,0){1.31}}
\put(5.4,1.24){\footnotesize$\emptyset$}

\put(6.82,0.4){$\frac{31}{63}$}
\put(6.37,-0.2){\footnotesize $\approx 0.492$}
\put(6.04,-0.8){\footnotesize $\psi'(\mathcal C_5\rtimes_\iota\mathcal C_2)$}
\put(7,1){\circle*{0.15}}

\put(7,1.35){\line(1,0){0.81}}
\put(7,1.35){\color{gray}\circle*{0.15}}
\put(8.9,1.35){\line(-1,0){0.75}}
\put(7.9,1.24){\footnotesize$\emptyset$}

\put(8.87,0.4){$\frac{1}{2}$}
\put(8.9,0.92){\begin{tikzpicture}
      \node[draw=black, fill=white, circle, inner sep=1.5pt] at (0,0) {};
     \end{tikzpicture}}
\put(9.14,1){\circle*{0.1}}
\put(9.34,1){\circle*{0.1}}
\put(9.65,1){\circle*{0.1}}
\put(10.5,1){\circle*{0.15}}
\put(10.32,0.4){$\frac{23}{43}$}
\put(8.47,-0.2){\footnotesize $\approx0.5\leftarrow\;\approx 0.534$}
\put(8.68,-0.8){\footnotesize $\psi'(\mathcal{C}_{2^{k-1}}\!\times \mathcal{C}_2)$}
\put(8.89,-1.4){\footnotesize $\psi'\big(M(2^k)\big)$}

\put(8.9,1.265){\begin{tikzpicture}
      \node[draw=black, fill=white, circle, inner sep=1.5pt] at (0,0) {};
     \end{tikzpicture}}
     
\put(10.5,1.35){\line(1,0){0.71}}
\put(10.5,1.35){\line(-1,0){1.45}}
\put(9.14,1.35){\color{gray}\circle*{0.1}}
\put(9.34,1.35){\color{gray}\circle*{0.1}}
\put(9.65,1.35){\color{gray}\circle*{0.1}}
\put(10.5,1.35){\color{gray}\circle*{0.15}}

\put(11.3,1.24){\footnotesize$\emptyset$}
\put(12.2,1.35){\line(-1,0){0.65}}
\put(12.2,1.265){\begin{tikzpicture}
      \node[draw=black, fill=white, circle, inner sep=1.5pt] at (0,0) {};
     \end{tikzpicture}}
\put(12.35,1.35){\line(1,0){1.7}}

\put(12.44,1.35){\color{gray}\circle*{0.1}}
\put(12.64,1.35){\color{gray}\circle*{0.1}}
\put(12.95,1.35){\color{gray}\circle*{0.1}}
\put(14,1.35){\color{gray}\circle*{0.15}}


\put(12.16,0.4){$\frac{4}{7}$}
\put(12.2,0.92){\begin{tikzpicture}
       \node[draw=black, fill=white, circle, inner sep=1.5pt] at (0,0) {};
     \end{tikzpicture}}
\put(12.44,1){\circle*{0.1}}
\put(12.64,1){\circle*{0.1}}
\put(12.95,1){\circle*{0.1}}
\put(11.66,-0.2){\footnotesize $\approx0.571\leftarrow\; \approx0.619$}
\put(12.17,-0.8){\footnotesize $\psi'(\mathcal{C}_3\rtimes_\iota\mathcal{C}_{2^k})$}

\put(13.82,0.4){$\frac{13}{21}$}
\put(14,1){\circle*{0.15}}

\put(2,1.9){\line(1,0){2.5}}   
\put(4.66,1.9){\dashline{0.2}(0,0)(0.9,0)}
\put(2,1.9){\color{gray}\circle*{0.15}}
\put(2.15,2.07){\footnotesize $M$-groups}

\put(0,2.7){\line(1,0){2.5}}   
\put(2.66,2.7){\dashline{0.2}(0,0)(0.9,0)}
\put(0,2.7){\color{gray}\circle*{0.15}}
\put(0.15,2.87){\footnotesize Supersoluble Groups}

\end{picture}

\end{center}

\vspace{1.6cm}

Finally, what can be said about the groups whose $\psi'$-value lies in the interval $\big(\frac{31}{77},\frac{19}{43}\big)$? It turns out that our techniques are far more general than they might initially appear. 
Indeed, once the classes of groups in the range $\big[\frac{19}{43},1\big]$ have {\it somehow} been identified, the proofs of~Theorems~A and B may, in principle, be carried out at the same time by means of an inductive argument essentially independent of the considerations on $M$-groups that we have made (details are provided in Remark~\ref{changeinduction}). However, we chose not to proceed in this way for two reasons: on the one hand, our aim was to introduce the class of $M$-groups into the study of the $\psi$-function for the first time; on the other hand, our approach enables us to explain the genesis of the groups that fall in the interval $\big[\frac{19}{43},1\big]$ (by exploiting the structure of finite $M$-groups and, remarkably, avoiding the use of algebra computational software or lists of groups of large order). Now, the same goes for the interval $\big(\frac{31}{77},1\big]$. In fact, having identified the missing classes of groups (that is, $Q_{16}$, $\mathcal C_3\times\mathcal C_3$, $\mathcal C_3\rtimes Q_8$ where the action is the unique non-trivial one, $D_{12}\simeq \mathcal{C}_6\rtimes_\iota \mathcal{C}_2$, $D_{18} \simeq \mathcal{C}_9 \rtimes_\iota \mathcal{C}_2$, $\mathcal C_5\rtimes_\iota C_{2^n}$, and their direct products with coprime cyclic groups) one can then use our arguments (dealing with several more similar cases) to show that the groups with $\psi'$-value in $\big(\frac{31}{77},1\big]$ are precisely those we just mentioned, together with those described in Theorems A and~B (a more detailed explanation is provided at the end of the paper). This gives a complete description of the groups satisfying the supersolubility criterion of~\cite{SSCriterion}. Moreover, we mention that it was established that, also in~\cite{SSCriterion}, the \mbox{non-supersoluble} groups with \hbox{$\psi'$-value} equal to~$\frac{31}{77}$ are precisely $A_4$ and its direct products with coprime cyclic groups.
As a concluding comment, with regard to our modularity criterion, it is interesting to notice that the additional~\hbox{$M$-groups} in the interval $\big(\frac{31}{77}$,$\frac{19}{43}\big]$ are: $\mathcal C_3\times\mathcal C_3$, $\mathcal C_5\rtimes_\iota C_{2^n}$, $D_{14}\simeq \mathcal{C}_7\rtimes_\iota \mathcal{C}_2$, and their direct products with coprime cyclic groups. Hence, we can affirm that the majority of the classes of groups identified by the supersolubility criterion of \cite{SSCriterion} are actually \hbox{$M$-groups,} a fact that further justifies our focus on $M$-groups.

\vspace{-1.2cm}

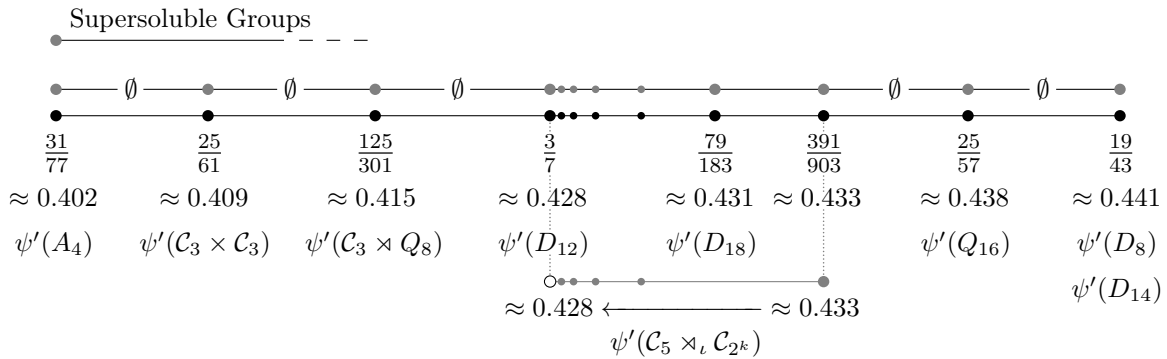
\begin{center}
\setlength{\unitlength}{1cm}
\begin{picture}(14,4)

\put(0,1){\line(1,0){14}}


\put(-0.18,0.4){$\frac{31}{77}$}
\put(-0.62,-0.2){\footnotesize $\approx 0.402$}
\put(-0.54,-0.8){\footnotesize $\psi'(A_4)$}
\put(0,1){\circle*{0.15}}

\put(0,1.35){\line(1,0){0.81}}
\put(0.9,1.24){\footnotesize$\emptyset$}
\put(2,1.35){\line(-1,0){0.85}}
\put(0,1.35){\color{gray}\circle*{0.15}}


\put(1.83,0.4){$\frac{25}{61}$}
\put(1.36,-0.2){\footnotesize $\approx 0.409$}
\put(2,1){\circle*{0.15}}
\put(1.1,-0.8){\footnotesize $\psi'(\mathcal C_3\times\mathcal C_3)$}

\put(2,1.35){\line(1,0){0.91}}
\put(4.2,1.35){\line(-1,0){0.94}}
\put(3,1.24){\footnotesize$\emptyset$}
\put(2,1.35){\color{gray}\circle*{0.15}}

\put(3.94,0.4){$\frac{125}{301}$}
\put(3.57,-0.2){\footnotesize $\approx 0.415$}
\put(3.28,-0.8){\footnotesize $\psi'(\mathcal C_3\rtimes Q_8)$}
\put(4.2,1){\circle*{0.15}}

\put(4.2,1.35){\line(1,0){0.91}}
\put(6.5,1.35){\line(-1,0){1.05}}
\put(5.2,1.24){\footnotesize$\emptyset$}
\put(4.2,1.35){\color{gray}\circle*{0.15}}

  \put(6.5,1){%
    \multiput(0,0)(0,-0.05){5}{\color{gray}\circle*{0.021}}
  }
  \put(6.5,0.2){%
    \multiput(0,0)(0,-0.05){5}{\color{gray}\circle*{0.021}}
  }
  \put(6.5,-0.25){%
    \multiput(0,0)(0,-0.05){6}{\color{gray}\circle*{0.021}}
  }
  \put(6.5,-0.85){%
    \multiput(0,0)(0,-0.05){6}{\color{gray}\circle*{0.021}}
  }

  \put(10.1,1){%
    \multiput(0,0)(0,-0.05){5}{\color{gray}\circle*{0.02}}
  }
    \put(10.1,0.2){%
    \multiput(0,0)(0,-0.05){5}{\color{gray}\circle*{0.02}}
  }
      \put(10.1,-0.25){%
    \multiput(0,0)(0,-0.05){20}{\color{gray}\circle*{0.02}}
  }

\put(6.5,-1.2){\color{gray}\line(1,0){3.66}}
\put(10.1,-1.2){\color{gray}\circle*{0.15}}
\put(6.42,-1.28){\begin{tikzpicture}
      \node[draw=black, fill=white, circle, inner sep=1.5pt] at (0,0) {};
    \end{tikzpicture}}

\put(5.9,-1.65){\footnotesize $\approx0.428  \xleftarrow{\hspace{1.92cm}} \;\approx 0.433$}

\put(7.29,-2.1){\footnotesize $\psi'(\mathcal C_5\rtimes_\iota\mathcal C_{2^{k}})$}

\put(6.65,-1.2){\color{gray}\circle*{0.1}}
\put(6.81,-1.2){\color{gray}\circle*{0.1}}
\put(7.1,-1.2){\color{gray}\circle*{0.1}}
\put(7.7,-1.2){\color{gray}\circle*{0.1}}

\put(6.38,0.4){$\frac{3}{7}$}
\put(6.5,1){\circle*{0.15}}
\put(5.82,-0.2){\footnotesize $\approx0.428$}

\put(5.82,-0.8){\footnotesize $\psi'(D_{12})$}

\put(6.65,1){\circle*{0.1}}
\put(6.81,1){\circle*{0.1}}
\put(7.1,1){\circle*{0.1}}
\put(7.7,1){\circle*{0.1}}
\put(8.67,1){\circle*{0.15}}
\put(8.41,0.4){$\frac{79}{183}$}
\put(8.02,-0.2){\footnotesize $\approx0.431$}
\put(8.03,-0.8){\footnotesize $\psi'(D_{18})$}
\put(10.1,1){\circle*{0.15}}
\put(9.84,0.4){$\frac{391}{903}$}
\put(9.43,-0.2){\footnotesize $\approx0.433$}



\put(10.1,1.35){\line(1,0){0.78}}
\put(12,1.35){\line(-1,0){0.78}}
\put(10.95,1.24){\footnotesize$\emptyset$}
\put(12,1.35){\line(1,0){0.81}}
\put(12.9,1.24){\footnotesize$\emptyset$}
\put(12,1.35){\color{gray}\circle*{0.15}}

\put(14,1.35){\line(-1,0){0.84}}
\put(14,1.35){\color{gray}\circle*{0.15}}


\put(6.5,1.35){\line(1,0){3.6}}
\put(6.5,1.35){\color{gray}\circle*{0.15}}
\put(6.65,1.35){\color{gray}\circle*{0.1}}
\put(6.81,1.35){\color{gray}\circle*{0.1}}
\put(7.1,1.35){\color{gray}\circle*{0.1}}
\put(7.7,1.35){\color{gray}\circle*{0.1}}
\put(8.67,1.35){\color{gray}\circle*{0.15}}
\put(10.1,1.35){\color{gray}\circle*{0.15}}

\put(11.82,0.4){$\frac{25}{57}$}
\put(12,1){\circle*{0.15}}
\put(11.37,-0.2){\footnotesize $\approx 0.438$}
\put(11.38,-0.8){\footnotesize $\psi'(Q_{16})$}

\put(13.82,0.4){$\frac{19}{43}$}
\put(14,1){\circle*{0.15}}
\put(13.37,-0.2){\footnotesize $\approx 0.441$}
\put(13.43,-0.8){\footnotesize $\psi'(D_8)$}
\put(13.36,-1.4){\footnotesize $\psi'(D_{14})$}


\put(0,2){\line(1,0){3}}
\put(3.2,2){\dashline{0.2}(0,0)(0.9,0)}
\put(0,2){\color{gray}\circle*{0.15}}
\put(0.18,2.17){\footnotesize Supersoluble Groups}

\end{picture}

\end{center}

\vspace{3cm}

\section{Preliminaries}
The aim of this section is to highlight, for the reader's convenience, results and notation that will be used frequently in our discussion, mainly concerning the function $\psi$ and groups whose subgroup lattice is modular. Some of the statements are rephrasings of known results, specifically crafted to emphasize useful consequences that we need.
\vspace{0.2cm}\\
\noindent We begin by listing the notation employed for the groups that arise in our analysis (here $n$ and $k$ denote integers):

\begin{itemize}[leftmargin=15pt]
    \item $\mathcal{C}_n$ denotes the cyclic group of order $n$ ($n\geq 1$).
    \item $D_{2n}$ denotes the dihedral group of order $2n$ ($n\geq 3$).
    \item If $H$ is an abelian group and $K$ is a cyclic group of even order, the notation 
    $$H\rtimes_\iota K$$
    denotes the semidirect product of $H$ by $K$, where the action $\iota: K \rightarrow \operatorname{Aut(H)}$ is defined by mapping a generator of $K$ to the inversion automorphism of $H$ (the choice of generator being immaterial). Thus, for example, $D_{2n} \simeq \mathcal C_{n}\rtimes_\iota\mathcal{C}_2$. With the exception of $D_8$, we shall frequently denote dihedral groups in this way, as this is the natural form in which they arise in our arguments; in particular, in the statements of~The\-o\-rems~A and B we find the groups $\mathcal{C}_5 \rtimes_\iota \mathcal{C}_2 \simeq D_{10}$ and $\mathcal{C}_7 \rtimes_\iota \mathcal{C}_2 \simeq D_{14}$. Moreover, it is worth noting that in the statement of Theorem A, we also find the group $\mathcal{C}_5 \rtimes_\iota \mathcal{C}_4$, which is more commonly known as the dicyclic group $\operatorname{Dic}_{20}$ of order~$20$.
    \item $Q_{2^k}$ denotes the generalized quaternion group of order $2^k$ ($k\geq 3$), i.e. the group defined by the presentation
    $$Q_{2^k} \coloneqq \langle x,y\mid x^{2^{k-2}}=y^2, \, y^4 = 1, \, yxy^{-1}=x^{-1} \rangle \; .$$
    \item $S_{2^k}$ denotes the semi-dihedral group of order $2^k$ ($k\geq 4$), i.e. the group defined by the presentation $$S_{2^k} \coloneqq \langle x,y\mid x^{2^{k-1}} = y^2 = 1, \, yxy^{-1}=x^{2^{k-2}-1} \rangle \; .$$
    \item $M(2^k)$ denotes the group of order $2^k$ ($k\geq 4$) defined by the presentation
    $$M(2^k) \coloneqq \langle x,y\mid x^{2^{k-1}} = y^2 = 1, \, yxy^{-1} = x^{2^{k-2} + 1} \rangle \; .$$ 
\end{itemize}   

In dealing with certain inequalities, the reader may encounter an annotation written above the inequality symbol $\leq$, such as
$$\overset{(\star)}{\leq},\; \overset{3.1}{\leq}\textnormal{ or } \overset{k=1}{\leq} \; ,$$
which indicates the justification for that step. Although the meaning of these expressions is likely self-explanatory, we illustrate them here for the sake of clarity. For instance,~$(\star)$ and~3.1 mean that the inequality follows from Equation $(\star)$ and the result~3.1, respectively. The annotation $k=1$ is more subtle: it is used when the left-hand side of the inequality involves a decreasing function of $k\geq 1$, which is then bounded above by its value at $k=1$, explicitly written on the right-hand side of the inequality. These concise and natural notations should improve the readability of the paper.

\subsection{The $\psi$-function}        

\begin{sloppypar}
\noindent One of the most elementary properties of the functions $\psi$ and $\psi'$, which is readily verified, is their multiplicativity on direct products with coprime 
factors: \hbox{if $A$ and $B$ are} finite groups such that $(|A|, |B|) = 1$, then~\hbox{$\psi(A \times B) = \psi(A) \cdot \psi(B)$} and~\hbox{$\psi'(A \times B) = \psi'(A) \cdot \psi'(B)$.} This fact is essential for computing the $\psi$-values and the $\psi'$-values of many groups, so we shall generally use it without further notice. On the other hand, the $\psi$-function does not generally behave well with respect to arbitrary products and quotients, and the first results we recall will help us handle certain specific cases of these operations.
\end{sloppypar}

\begin{lemma}[\cite{GeneralizedBoundPsiCyclic}, Proposition 2.6] 
 \label{BoundPsiQuotient} 
    Let $G$ be a finite group and $N$ a normal subgroup of~$G$. Then $\psi(G) \leq \psi(G/N)\cdot|N|^2$.
\end{lemma}

\begin{lemma}\label{BoundPsiPNormalCyclicSylow} 
Let $G$ be a finite group, and let $P$ be a normal cyclic Sylow subgroup of $G$. The following results hold.
\begin{enumerate}[label=$(\roman*)$]
\item \emph{(\cite{FIRST}, Corollary B)} $\psi(G)\leq\psi(P)\cdot\psi(G/P)$, with equality if and only if $P$ is central in~$G$.
\label{BoundPsiPNormalCyclicSylowBound}
\item \emph{(\cite{NonCyclicCriterionI}, Lemma 2.2(5))} If $G=P\rtimes H$ for some subgroup $H$ of $G$, then $$\psi(G) = |P|\cdot\psi(H) + \left(\psi(P)-|P|\right)\cdot\psi(C_H(P)) \; .$$
\label{BoundPsiPNormalCyclicSylowFormulaSDP}
\end{enumerate}
\end{lemma}
Next, we need some expressions and bounds for the value of $\psi$ in the case of a cyclic group. In particular, the first lemma provides an explicit formula for the $\psi$-value of a cyclic group in terms of the parameters of the prime factorization of its order, which follows easily from the multiplicativity of $\psi$ together with an obvious induction argument.

\begin{lemma}\label{psiCyclic}
 Let $n=p_1^{\alpha_1}\cdot\ldots\cdot p_k^{\alpha_k}$ with $k,\alpha_1,\ldots,\alpha_k$ positive integers and $p_1,\ldots,p_k$ distinct prime numbers. Then
 $$\psi(\mathcal{C}_n) = \frac{p_1^{2\alpha_1+1}+1}{p_1+1}\cdot \ldots \cdot \frac{p_k^{2\alpha_k+1}+1}{p_k+1} \; .$$
\end{lemma}

\begin{lemma}[\cite{GeneralizedBoundPsiCyclic}, Lemma 2.3] \label{BoundPsiCn}
    Let $n=p_1^{\alpha_1}\cdot\ldots\cdot p_k^{\alpha_k}$, where $k,\alpha_1, \ldots, \alpha_k$ are positive integers and $p_i$ are prime numbers such that $p_1 < \ldots <  p_k$. Let $q_1,  q_2, \ldots, q_k$ be the first $k$ prime numbers. The following hold.
    \begin{enumerate}[label=$(\roman*)$]
    \item \label{BoundPsiCn_ONE} We have $$\psi(\mathcal{C}_n) > \frac{p_1}{p_1+1}\cdot \ldots \cdot \frac{p_k}{p_k+1}\cdot n^2\geq \frac{q_1}{q_1+1}\cdot \ldots \cdot \frac{q_k}{q_k+1}\cdot  n^2 \; .$$
    \item \label{BoundPsiCn_TWO} If $t\in\{1,\ldots, k\}$, then
     $$\psi(\mathcal{C}_n) > \frac{p_1}{p_1+1}\cdot \ldots \cdot \frac{p_{t-1}}{p_{t-1}+1} \cdot \frac{p_t}{p_k+1}\cdot n^2\geq \frac{q_1}{q_1+1}\cdot \ldots \cdot \frac{q_{t-1}}{q_{t-1}+1} \cdot \frac{q_t}{p_k+1}\cdot  n^2 \; .$$
In particular,
$$\psi(\mathcal{C}_n) > \frac{p_1}{p_k+1}\cdot n^2 \; ,$$
where $p_1$ and $p_k$ are, respectively, the smallest and the largest prime dividing $n$.
\end{enumerate}
\end{lemma}            

Here, as in many previous papers, we are interested in obtaining information about groups whose $\psi'$-value exceeds certain critical thresholds. When such a threshold is sufficiently large, the group's structure is strongly constrained by the existence of a very large cyclic subgroup. This is a direct consequence of the lower bounds provided by Lemma~\ref{BoundPsiCn}, which in turn yield corresponding upper bounds on the indices of a particular cyclic subgroup; the argument is standard and proceeds as follows. Let $G$ be a finite group and suppose that \mbox{$\psi'(G)>$ \scalebox{0.8}{$\dfrac{r}{s}$}} for some positive integers $r$ and $s$. From this inequality, using a bound of the form $\psi(\mathcal{C}_{|G|}) > k \cdot |G|^2$ for a suitable real number $k$, such as those provided by Lemma~\ref{BoundPsiCn}, one obtains \mbox{\scalebox{0.8}{$\dfrac{\psi(G)}{|G|}$} $> k \cdot |G|$.} Having recognized that the left-hand side is just the arithmetic mean of the orders of the elements of $G$, the generalized pigeonhole principle ensures the existence of an element whose order exceeds the right-hand side, thus leading to the next lemma (see, for instance, Lemma 2.1 in \cite{SolvableCriterion}).

\begin{lemma}\label{CyclBoundIndex}
Let $G$ be a finite group of order $n=p_1^{\alpha_1}\ldots p_k^{\alpha_k}$, where $k,\alpha_1, \ldots, \alpha_k$ are positive integers and $p_i$ are prime numbers such that $p_1 < \ldots <  p_k$. Let $q_1, q_2, \ldots, q_k$ be the first $k$ prime numbers. If $\psi'(G) > \dfrac{r}{s}$, then the following hold. 
\begin{enumerate}[label=$(\roman*)$]
\item \label{CyclBoundIndex_ONE} There exists $x\in G$ such that $$|G:\langle x \rangle| < \frac{s}{r}\cdot\frac{p_1+1}{p_1}\cdot \ldots  \cdot \frac{p_k+1}{p_k}\leq \frac{s}{r}\cdot\frac{q_1+1}{q_1}\cdot \ldots \cdot \frac{q_k+1}{q_k} \; .$$
\item \label{CyclBoundIndex_TWO} If $t\in\{1,\ldots, k\}$, then there exists $x\in G$ such that $$|G:\langle x \rangle| < \frac{s}{r}\cdot\frac{p_1+1}{p_1}\cdot \ldots \cdot \frac{p_{t-1}+1}{p_{t-1}} \cdot \frac{p_k+1}{p_t}\leq \frac{s}{r}\cdot\frac{q_1+1}{q_1}\cdot \ldots \cdot \frac{q_{t-1}+1}{q_{t-1}} \cdot \frac{p_k+1}{q_t} \; .$$
 \end{enumerate}
\end{lemma}

\begin{sloppypar} In fact, Lemma \ref{CyclBoundIndex} will often be used in combination with the following result \mbox{of~Luc\-chi\-ni.}
\end{sloppypar}

\begin{lemma}[\cite{Lucchini}, Theorem 2.20] \label{lucchini}
Let $A$ be a proper cyclic subgroup of a finite group~$G$. If $K=\mathrm{core}_G(A)$, then $|A:K|<|G:A|$.
\end{lemma}

Finally, the best possible upper bound for the $\psi'$-value among non-cyclic groups is given by the following fundamental result. 

\begin{theorem}[\cite{NonCyclicCriterionII}, Theorem 4] \label{BoundPsiNonCycl}
Let $G$ be a non-cyclic group and let $q$ be the smallest prime dividing the order of $G$. Then
$$\psi'(G) \leq f(q) \coloneqq \frac{\left[(q^2-1)q+1\right](q+1)}{q^5+1} < \frac{1}{q-1} \; ,$$
\begin{sloppypar}\noindent with equality if and only if \hbox{$G \simeq \left( \mathcal{C}_q \times \mathcal{C}_q \right) \times \mathcal{C}_m$} for some positive integer $m$ such that~\hbox{$(m,q!)=1$.} In particular, since the function $f$ is strictly decreasing on $[2, +\infty)$, the maximum value of the function $\psi'$ among non-cyclic groups is $f(2)=\frac{7}{11}$, attained only when $q=2$ and~\hbox{$G \simeq \left( \mathcal{C}_2 \times \mathcal{C}_2 \right) \times \mathcal{C}_m$} for some odd positive integer $m$.
\end{sloppypar}
\end{theorem}

\subsection{$M$-groups}

In this subsection, we recall some properties and results concerning $M$-groups. Before delving into specific theorems, let us first recall the definition. A group $G$ is said to be an \emph{$M$-group} if its subgroup lattice is modular. This means that for all subgroups $ H, K, L \leq G $ with $ H \leq L $, the equality
$$
\langle H,\,K \cap L \rangle = \langle H,\,K \rangle \cap L
$$
holds. This identity, commonly known as {\it Dedekind’s modular law}, expresses a weak form of distributivity of the subgroup lattice and always holds when $K$ is normal, so in particular every abelian group is an $M$-group. Thus, $M$-groups can be viewed as a natural generalization of Dedekind groups in which normality of all subgroups is replaced by modularity. In this sense, $M$-groups provide the next level of complexity beyond~De\-de\-kind groups, while preserving enough order to permit a detailed structural theory. For a comprehensive treatment of $M$-groups, see \cite{Schmidt}; here, we just present some essential information. To begin with, we mention that the class of $M$-groups is closed with respect to forming subgroups, quotients and direct products of finite coprime $M$-groups. We are also going to state the complete structure theorem for finite $M$-groups, but in order to do so, we first need to recall two relevant related concepts: Iwasawa triples and $P^*$-groups.
\begin{itemize}[leftmargin=15pt]
    \item Let $G$ be a finite $p$-group, where $p$ is a prime number. An {\it Iwasawa triple} for $G$ is a triple~$\big(A,\,b, \, s \big)$, where $A$ is a normal abelian subgroup of $G$, $b$ is an element of $G$ and~$s$~is a positive integer, with $s\geq 2$ in case $p=2$, such that $G=A\langle b \rangle$ and $a^b = a^{1+p^s}$ for all~$a\in A$. If such a triple exists for $G$, then we say that $G$ admits an Iwasawa triple. Groups admitting an Iwasawa triple are $M$-groups (\cite{Schmidt}, as a consequence of~Lem\-ma~2.3.4).
    
    \item Let $p$ and $q$ be prime numbers such that $q$ divides $(p-1)$ (so $q<p$), and let $n,k$ be positive integers. We define the group $P^*(p^n,q^k)$ as the semidirect product of an elementary abelian group $A$ of order $p^n$ by a cyclic group $\langle x\rangle$ of order $q^k$, where the action of~$x$ on~$A$ is an automorphism of $A$ of order $q$ that fixes the subgroups (a \mbox{so-called} {\it power automorphism}) --- clearly, such automorphisms exist, and any two of them produce the same group up to isomorphism. A finite group is said to be a {\it $P^*$-group} if and only if it is isomorphic to one of these groups. Detailed information on~$P^\ast$-groups can be found in \cite{Schmidt}, \textsection 2.4. Here, we only need to mention that $P^*$-groups are $M$-groups (\cite{Schmidt}, Lemma 2.4.1), and to observe that if $q = 2$, then the action of $x$ on $A$ is given by the inversion automorphism of~$A$, that is, in our notation, $P^*(p^n, 2^k) = A \rtimes_\iota \langle x \rangle$ (since power automorphisms of a finite abelian group map all elements to the same power).        
\end{itemize}
\noindent Now, we can state the complete structure theorem for finite $M$-groups.

\begin{theorem}[\cite{Schmidt}, Theorems 2.3.1 and 2.4.4]\label{StructureFiniteMgroups}
A finite group is an M-group if and only if it is isomorphic to a direct product, with pairwise coprime factors, of finitely many groups of the following types\textnormal:\\

$\left.\begin{array}{l}
\hspace{-0.7cm}\bullet \ \  Q_8\times A, \text{where $A$ is a finite elementary abelian $2$-group\textnormal;}\\[8pt]
\hspace{-0.7cm}\bullet \ \ \text{groups admitting an Iwasawa triple\textnormal;}
\end{array}\right\}$ \textnormal(finite primary $M$-groups\textnormal)\\[3pt]

\hspace{-0.47cm}$\bullet \ \ P^\ast$-groups.
\end{theorem}   

In the definition of a $P^*$-group given above, all subgroups of $A$ are readily seen to be normal, and as a consequence every $P^*$-group is supersoluble. Finite primary \hbox{$M$-groups} are obviously supersoluble as well, and so Theorem~\ref{StructureFiniteMgroups} implies that all finite $M$-groups are supersoluble. On the other hand, the dihedral group $D_8$ of order $8$ shows that supersoluble groups need not be $M$-groups. Indeed, writing $D_8 = \langle a \rangle \rtimes_\iota \langle  x \rangle$, with $\langle a \rangle \simeq \mathcal{C}_4$ and $\langle x \rangle \simeq \mathcal{C}_2$, the three subgroups $H=\{1,x\}, K=\{1,ax\}, L=\{1,a^2,x,a^2x\}$ do not satisfy Dedekind's modular law. Moreover, since the other non-abelian groups of order at most~$8$ are $S_3$ (which is an $M$-group because its subgroup lattice is isomorphic to that of the abelian group $\mathcal{C}_3\times \mathcal{C}_3$), and $Q_8$ (which is an $M$-group because all its subgroups are normal), we see that $D_8$ is the smallest group whose subgroup lattice is not modular. Hence, its $\psi'$-value is so essential in our work that, in closing this paragraph, we explicitly compute it. In the above notation, in $D_8$ we count the identity, five elements of order $2$ (namely, $a^2, x, ax, a^2x, ax^3$) and two elements of order $4$ (namely, $a, a^3$), thus obtaining~$\psi(D_8)=19$ and then, by using Lemma~\ref{psiCyclic}, \hbox{$\psi'(D_8)=$ \kern0.12em \scalebox{0.8}{$\dfrac{\psi(D_8)}{\psi({\mathcal{C}_8})}$}\kern0.3em $=$ \kern0.08em \scalebox{0.8}{$\dfrac{19}{43}$}.}

\section{Main results}

The aim of this section is to prove Theorems A and B, and show how the arguments employed in their proofs can be easily extended to deal with groups in the range $\big(\frac{31}{77},1\big]$. We start with some preliminary observations about the primes that always need to divide the order of the groups we are considering, and the description of the corresponding primary cases.

\begin{lemma}\label{smallestprime1943-3177}
Let $G$ be a finite non-cyclic group.
\begin{enumerate}[label=$(\roman*)$]
\item If $\psi'(G) \geq \frac{19}{43}$, then the order of $G$ is divided by $2$.
\label{smallestprime1943}
\item If $\psi'(G) > \frac{31}{77}$, then the order of $G$ is divided by either $2$ or $3$.
\label{smallestprime3177}
\end{enumerate}
\end{lemma}
\begin{proof}
    The statement follows directly from Theorem~\ref{BoundPsiNonCycl}, noting that \mbox{$f(3)=\frac{25}{61} < \frac{19}{43}$} and~\mbox{$f(5)=\frac{121}{521}<\frac{31}{77}$}.
\end{proof}

\begin{coroll}\label{Primary1943nc}
    A non-cyclic primary group with $\psi'$-value greater than or equal to $\frac{19}{43}$ is a $2$-group.
\end{coroll}

The following result gives a complete description of the finite $2$-groups with a~\hbox{$\psi'$-value} strictly greater than $\frac{31}{77}$. In particular, it establishes Theorem A for $2$-groups, and consequently, by Corollary~\ref{Primary1943nc}, for primary groups.

\begin{prop}\label{2g-classification}
Let $G$ be a finite $2$-group.

\begin{enumerate}[label=$(\roman*)$]
\item If $\psi'(G) > \frac{19}{43}$, then~$G$ is an $M$-group. More precisely, up to isomorphism, $G$ belongs to one of the following pairwise disjoint classes of groups\textnormal:
\begin{itemize}
    \item cyclic groups $\mathcal{C}_{2^k}$, where $k\geq 0$;
    \item $\mathcal{C}_{2^{k-1}} \!\times \mathcal{C}_2$, where $k\geq 2$;
    \item $M(2^k) \coloneqq \langle x,y \mid x^{2^{k-1}} = y^2 = 1, \, yxy^{-1} = x^{2^{k-2} +1} \rangle$, where $k\geq 4$;
    \item $Q_8$.
\end{itemize}
\label{2g>1943}
\item If $\psi'(G) = \frac{19}{43}$, then~$G\simeq D_8$.
\label{2g=1943}
\item If $\frac{31}{77} < \psi'(G) < \frac{19}{43}$, then~$G\simeq Q_{16}$.
\label{2g3177-1943}
    \end{enumerate}
\end{prop}

\begin{proof}

Set $|G|=2^k$, where $k\in\mathbb{N}_0$; clearly, we may assume $k>0$. 
Observe that under any of the assumptions in points (i), (ii), or (iii), we always have \hbox{$\psi'(G) > \frac{31}{77}$}, hence we adopt this as our working assumption.
By Lemma~\ref{CyclBoundIndex}\ref{CyclBoundIndex_ONE}, there exists $x\in G$ such that
    $$|G:\langle x \rangle| < \frac{77}{31}\cdot\frac{3}{2} \approx 3.72\ldots\ .$$
    
    If $G$ is cyclic, then $G \simeq \mathcal{C}_{2^k}$, which is an $M$-group with $\psi'$-value equal to~$1$. Assume otherwise, so $|G:\langle x \rangle|=2$; in other words, this means that $G$ is a finite non-cyclic~\hbox{$2$-group} with a maximal cyclic subgroup of index $2$. Therefore, up to isomorphism,~$G$ belongs to one of the following classes of groups (\cite{SuzukiBooks}, Vol. II, Chapter~4, Theorem 4.1): 
    \begin{itemize}
    \item $\mathcal{C}_{2^{k-1}} \!\times \mathcal{C}_2$, where $k \geq 2$;
    \item $M(2^k) \coloneqq \langle x,y\mid x^{2^{k-1}} = y^2 = 1, \, yxy^{-1} = x^{2^{k-2} + 1} \rangle$, where $k\geq 4$;
    \item $D_{2^k} \coloneqq \langle x,y\mid x^{2^{k-1}} = y^2 = 1, \, yxy^{-1}=x^{-1}\rangle$, where $k\geq 3$;
    \item $Q_{2^k} \coloneqq \langle x,y\mid x^{2^{k-2}}=y^2, \, y^4 = 1, \, yxy^{-1}=x^{-1} \rangle$, where $k\geq 3$;
    \item $S_{2^k} \coloneqq \langle x,y\mid x^{2^{k-1}} = y^2 = 1, \, yxy^{-1}=x^{2^{k-2}-1} \rangle$, where $k\geq 4$.
    \end{itemize}
    
    Now, in order to establish point (i) of the statement, we proceed to determine which of these groups have a $\psi'$-value strictly greater than $\frac{19}{43}$. Then, whenever such a group is identified, we verify that it is an $M$-group, as predicted. 
    
    First, the computation
    $$\psi(\mathcal{C}_{2^{k-1}} \!\times \mathcal{C}_2) = \frac{2^{2k}+5}{3}$$ follows straightforwardly from the observation that, given $a\in \mathcal{C}_{2^{k-1}}\setminus\{1\}$ and $b\in \mathcal{C}_2$, we have $o(ab)=o(a)$. Hence, we observe that the inequality
 $$\psi'(\mathcal{C}_{2^{k-1}} \!\times \mathcal{C}_2) = \frac{\psi(\mathcal{C}_{2^{k-1}} \!\times \mathcal{C}_2)}{\psi(\mathcal{C}_k)} = \frac{2^{2k}+5}{2^{2k+1}+1} > \frac{19}{43}$$
holds for every $k\geq 2$. Clearly, the groups $\mathcal{C}_{2^{k-1}}\!\times\mathcal{C}_2$ with $k\geq 2$ are $M$-groups, since they are abelian.

 Next, the computation of the $\psi$-value in the other cases, where an explicit presentation in terms of two generators $x$ and $y$ is provided, can be carried out by observing that half of the $2^k$ elements are those of $\langle x \rangle$ (allowing the use of Lem\-ma~\ref{psiCyclic}), while the other half are the elements $x^t y$ with $0\leq t \leq 2^{k-1}-1$ (whose orders can be computed directly from the presentation). Thus, we obtain:

\begin{itemize}
    \item $\displaystyle \psi'(M(2^k)) = \frac{2^{2k}+5}{2^{2k+1}+1}>\frac{19}{43}\, ,$\vspace{0.15cm} which is true for every $k\geq 4$. These are~\hbox{$M$-groups} because they are $2$-groups admitting the Iwasawa triple $\Big(\langle x \rangle,\, y, \, k-2 \Big)$.
    \item $\displaystyle \psi'(D_{2^k}) = \frac{2^{2k-1}+3\cdot 2^k + 1}{2^{2k+1}+1} > \frac{19}{43} \, ,$ which is false for every $k\geq 3$.
    \item $\displaystyle \psi'(Q_{2^k}) = \frac{2^{2k-1}+3\cdot 2^{k+1} +1}{2^{2k+1}+1} > \frac{19}{43} \, ,$\vspace{0.15cm} which is false for every $k\geq 4$ and true for~\hbox{$k=3$}, which corresponds to $Q_8$. Clearly, $Q_8$ is an $M$-group, as all of its subgroups are normal.
    \item $\displaystyle \psi'(S_{2^k}) = \frac{2^{2k-1}+9\cdot 2^{k-1} +1}{2^{2k+1}+1} > \frac{19}{43} \, ,$ which is false for every $k\geq 4$.
\end{itemize}

Point (i) is thus proved. Points (ii) and (iii) follow analogously, by imposing the \hbox{$\psi'$-values} of the previously mentioned classes of groups equal to $\frac{19}{43}$ and strictly greater than $\frac{31}{77}$, respectively.
\end{proof}

Next, we present an analogue of Proposition \ref{2g-classification} for $3$-groups.

\begin{prop}\label{3g-classification}
Let $G$ be a finite $3$-group with $\psi'(G) > \frac{31}{77}$. Then either $G\simeq\mathcal{C}_{3^t}$ for some integer $t \geq 0$, or $G\simeq \mathcal C_3\times\mathcal C_3$.
\end{prop}
\begin{proof}
Set $|G|=3^t$, where $t$ is a positive integer. By Lemma \ref{CyclBoundIndex}\ref{CyclBoundIndex_ONE}, there exists $x\in G$ such that
    $$|G:\langle x \rangle| < \frac{77}{31}\cdot\frac{4}{3} \approx 3.31\ldots\ .$$
Clearly, we may assume that $G$ is non-cyclic, so $|G:\langle x \rangle|=3$. 

If $G$ is abelian, then $G\simeq\mathcal{C}_{3^{t-1}} \!\times\mathcal{C}_3$. Since for all $a\in\mathcal{C}_{3^{t-1}}\setminus\{1\}$ and~$b\in\mathcal{C}_3$, we have~$o(ab)=o(a)$, it easily follows that 
$$\psi'(\mathcal{C}_{3^{t-1}}\!\times \mathcal{C}_3)=\frac{3^{2t}+19}{3^{2t+1}+1} \; ,$$
which is strictly greater than $\frac{31}{77}$ if and only if $t=2$; in this case, $G\simeq \mathcal{C}_3 \times \mathcal{C}_3$.

Now assume that $G$ is not abelian. Then, by Theorem~4.1 of~\cite{SuzukiBooks}, Vol.~II, Chap\-ter~4,~$G$ is isomorphic to the group
$$M(3^t)\coloneqq \langle x,y \, \mid \, x^{3^{t-1}} = y^3 = 1, \, yxy^{-1} = x^{1+3^{t-2}} \rangle = \langle x \rangle \rtimes \langle y \rangle \, ,$$
where $t\geq 3$.
By direct computation, using the given presentation, it is straightforward to verify that for all $a\in\langle x\rangle \setminus\{1\}$ and $b\in \langle y \rangle$, we have $o(ab)=o(a)$, in analogy with the behavior observed in the group $\mathcal{C}_{3^{t-1}} \!\times \mathcal{C}_3$. Consequently, computing the $\psi$-value once again yields
$$\psi'(M(3^t)) = \frac{3^{2t}+19}{3^{2t+1}+1} \; ,$$
which, as before, leads to $t=2$, now a contradiction. The statement is proved.
\end{proof}

\begin{remark}
In proving Propositions \ref{2g-classification} and \ref{3g-classification}, one may observe that, for each $n\geq 2$, we have $\psi(\mathcal{C}_{2^{n-1}}\!\times \mathcal{C}_2) = \psi(M(2^n))$ and $\psi(\mathcal{C}_{3^{n-1}}\!\times \mathcal{C}_3) = \psi(M(3^n))$. Indeed, these are special cases of a more general result. For a prime number $p$ and a positive integer $n\geq 3$, consider the group defined by the following presentation:
$$M(p^n) \coloneqq \langle x,y \, \mid \,  x^{p^{n-1}}=y^p=1, \, yxy^{-1}=x^{1+p^{n-2}} \rangle \; .$$
Clearly, $M(p^n) \simeq \mathcal{C}_{p^{n-1}} \rtimes_{\theta} \mathcal{C}_p$, where $\theta$ is the action defined by mapping an element to its $(1+p^{n-2})$-th power. 
We claim that
$$\psi(\mathcal{C}_{p^{n-1}}\!\times \mathcal{C}_p) = \psi(M(p^n)) \; .$$
For $p = 2$ and $p = 3$, this has already been established, essentially by direct computation. Now, we are going to treat the case of an arbitrary odd prime $p$, using a different argument. It is clear, without any computation, that the derived subgroup of $M(p^n)$ is~\mbox{$\langle [x,y] \rangle = \langle x^{p^{n-2}}\rangle$} and has order $p$, hence it is central. Consequently, for all $a,b\in M(p^n)$, we have
$$\displaystyle (ab)^p = a^p b^p [b,a]^{\binom{p}{2}} = a^p b^p \; .$$
In particular, in $M(p^n)$ the $p$-th powers distribute over products of elements in the same way as in $\mathcal{C}_{p^{n-1}}\!\times \mathcal{C}_p$, so the claim is immediate. This provides further examples of \mbox{non-isomorphic} groups with the same $\psi$-value.
\end{remark}

Now, we move to non-primary groups. \mbox{Proposition 2.1} in~\cite{NilpotentCriterion} states that any finite group that is not a $2$-group and whose $\psi'$-value exceeds $\frac{13}{21}$ has a normal cyclic~Sy\-low~\mbox{$r$-subgroup} for some $r \in \{2, p\}$, where $p$ is the largest prime dividing the group order; moreover, it can be shown that the case $r = p$ always occurs (see Corollary 1.2 in the same paper). Our next key result, namely Proposition~\ref{3177SylNormCycl}, is a generalization of that fact, obtained by relaxing the bound on the $\psi'$-value from $\frac{13}{21}$ to $\frac{31}{77}$. However, this broader generality comes at the cost of explicitly excluding the $3$-groups from the analysis, in addition to the $2$-groups already excluded; indeed, $\psi'(\mathcal{C}_3\times\mathcal{C}_3)= \frac{25}{61} >\frac{31}{77}$ (actually, $\mathcal{C}_3\times\mathcal{C}_3$ is the only $3$-group that needs to be excluded, see Proposition~\ref{3g-classification}).
Under the same conditions, this is the maximal generality that can be achieved: the restriction on the $\psi'$-value cannot be further lowered, as $\psi'(A_4)=\frac{31}{77}$ and~$A_4$ does not have Sylow subgroups that are both cyclic and normal.

A second point is worth discussing. Conjecture 1.1 in~\cite{SSCriterion}, proved therein, states the following supersolubility criterion: a finite group whose $\psi'$-value exceeds $\frac{31}{77}$ is supersoluble. Clearly, this implies that the group has a normal~Sy\-low~\hbox{$p$-sub}\-group, where $p$ is the largest prime dividing the group order. However, this Sylow subgroup need not be cyclic. To ensure cyclicity as well, Proposition~\ref{3177SylNormCycl} excludes~\hbox{$2$-groups} and $3$-groups in the hypothesis, while allowing for the alternative existence of a normal cyclic Sylow~\hbox{$2$-sub}\-group in the conclusion. More precisely, the groups for which these additional conditions are required are: $Q_8$, $D_8$, $Q_{16}$, $\mathcal{C}_2\times\mathcal{C}_2$, $\mathcal{C}_3\times \mathcal{C}_3$ and $(\mathcal{C}_3\times \mathcal{C}_3)\times \mathcal{C}_2$ (see the discussion at the end of the introduction). On the other hand, Proposition \ref{3177SylNormCycl} yields the supersolubility criterion via a straightforward induction argument. In view of these observations, Proposition \ref{3177SylNormCycl} may be regarded as a stronger variant of Conjecture~1.1 in~\cite{SSCriterion}. Indeed, to establish the former, we employ techniques stemming from the proof of the latter. As a preliminary step, we have an easy lemma that allows us to avoid both several computations with GAP~\cite{GAP} and the use of lists of groups of large order.

\begin{lemma}\label{techGAP}
Let $G$ be a finite group of order $n$ with a non-normal cyclic subgroup $X$ of order $m$. Moreover, let $K=\operatorname{core}_G(X)$.

\begin{enumerate}[label=$(\roman*)$]
    \item If $n=36$ and $m=12$, then~\mbox{$G\simeq \mathcal{C}_3\times(\mathcal{C}_3\rtimes_{\iota} \mathcal{C}_4)$} and $\psi(G)=243$.
    \label{techGAP1}
    \item If $n=54$ and $m=18$, then~\mbox{$G\simeq \mathcal{C}_9\times S_3$} and $\psi(G)=553$.
    \label{techGAP2}
    \item If $n=24$ and $m=6$, and moreover $|K|=2$, then one of the following holds\textnormal:
    \begin{itemize}[leftmargin=14pt]
    \item[\textbullet] \mbox{$G\simeq \mathcal{C}_2 \times A_4$} and $\psi(G)=87$\textnormal;
    \item[\textbullet] \mbox{$G\simeq Q_8 \rtimes \mathcal{C}_3$}, where the action is the unique \mbox{non-trivial} one, and $\psi(G)=99$.
    \end{itemize}
    \label{techGAP3}
    \vspace{-0.5\topsep}
    \item If $n=36$ and $m=9$, and moreover $|K|=3$, then~\mbox{$G\simeq (\mathcal C_2\times\mathcal C_2)\rtimes \mathcal{C}_9$} where the action is the unique \mbox{non-trivial} one, and~$\psi(G)=265$.
    \label{techGAP4}
\end{enumerate}
\end{lemma}
\begin{proof} First of all, the computations of the $\psi$-values of the groups under consideration are straightforward exercises. 
Nevertheless, even though this is not necessary, one can use GAP~\cite{GAP}: their ID's in the Small Group Library (\textit{SmallGrp}) are
$\mathtt{[36,6]}$, $\mathtt{[54,4]}$, $\mathtt{[24,13]}$, $\mathtt{[24,3]}$, and $\mathtt{[36,3]}$, respectively, in order of appearance.

\medskip

\noindent We start by proving (i). Since $|G:K|$ divides $|G:X|!=3!$, we deduce that $G/K \simeq S_3$ and $K\simeq\mathcal{C}_6$.
As $X\not\trianglelefteq G$, we have $X\neq C_G(K)$, and hence $X<C_G(K)$. Noting that $X$ is a maximal subgroup of~$G$, we obtain~$K\leq Z(G)$. Moreover, \mbox{$Z(G)/K\leq Z(G/K) = \{1\}$,} so~$K=Z(G)$.

Let~$H/K$ denote the subgroup of order $3$ in $G/K$, and let $P$ be a Sylow $3$-subgroup of~$H$. As~$|H:P|=2$, $P$ is characteristic in $H$, hence normal in $G$; in particular,~\mbox{$G=P\rtimes Q$,} where $Q\simeq\mathcal{C}_4$ is a Sylow $2$-subgroup of $G$. Then, given that $P$ is abelian (as it has order $9$), by Theorem~2.3 in Chapter 5 of \cite{Gorenstein}, we have
$$P=C_P(Q) \times [P,Q] = (P\cap Z(G)) \times [P,Q] \; .$$
Since $P\cap Z(G) = P \cap K$, it follows that $P\cap Z(G)\simeq \mathcal{C}_3$, upon examining the possible orders of $P\cap K$. Consequently, $[P,Q]\simeq\mathcal{C}_3$. Finally, since $[P,Q]\trianglelefteq G$, we obtain
$$G=\Big(\big(P\cap Z(G)\big) \times [P,Q]\Big) \rtimes Q = \big(P\cap Z(G)\big) \times \big([P,Q] \rtimes Q \big)\simeq\mathcal{C}_3\times (\mathcal{C}_3\rtimes_{\iota} \mathcal{C}_4) \; .$$

Next, we prove (ii). As in (i), we have $G/K\simeq S_3$ and \hbox{$K=Z(G)$}; clearly,~\mbox{$K\simeq\mathcal C_9$.} Let~$P$ be a Sylow $3$-subgroup of $G$. Then, $|P|=27$ and $|G:P|=2$, so~$P$ is normal in $G$; in particular, $G=P\rtimes Q$, where $Q\simeq\mathcal{C}_2$ is a Sylow $2$-subgroup of~$G$. Moreover, $K\leq P$ (since $P$ is the unique Sylow $3$-subgroup of $G$), so $P$ is \mbox{central-by-cyclic} and hence abelian.
Therefore, by Theorem~2.3 in Chapter~5 of~\cite{Gorenstein}, we obtain~\mbox{$P=C_P(Q)\times [P,Q] = K \times [P,Q]$,} and so~\mbox{$G = K \times ([P,Q] \rtimes Q) \simeq \mathcal{C}_9\times S_3$.}

\medskip

We now turn to the proof of (iii). Clearly, $G/K\simeq A_4$.
Since \mbox{$|G\!:\!C_G(K)|\!=\!\lvert\operatorname{Aut}(K)\rvert\!=\!1$,} we have $K\leq Z(G)$. Moreover, $Z(G)/K \leq Z(G/K)=Z(A_4)=\{1\}$, so $K=Z(G)$.

Let $P/K$ be the subgroup of order $4$ in $G/K$. Then $|P|=8$ and $P\trianglelefteq G$; in particular, $G=P\rtimes Q$, where $Q\simeq\mathcal{C}_3$ is a Sylow $3$-subgroup of $G$. If $Q\trianglelefteq G$, then $Q$ would be the unique subgroup of order $3$ of $G$, and so $X=KQ\trianglelefteq G$, a contradiction. Therefore, the action in $G=P\rtimes Q$ is non-trivial.

Assume that $P$ is not abelian. If $P\simeq D_8$, then $\lvert\operatorname{Aut}(D_8)\rvert = 8$ implies $G=P\times Q$, which is a contradiction. Thus $P\simeq Q_8$, and consequently $G\simeq Q_8\rtimes \mathcal{C}_3$, where the action is non-trivial (and the only possible one).

Now assume that $P$ is abelian. Then, Theorem~2.3 in Chapter~5 of~\cite{Gorenstein} yields
$$P= C_P(Q) \times [P,Q] = Z(G) \times [P,Q] \; ;$$
in particular,~$\big|[P,Q]\big|=4$. Consequently, $G= Z(G) \times \big( [P,Q] \rtimes Q\big)$.
If $[P,Q]\simeq \mathcal{C}_4$, then~\mbox{$[P,Q] \rtimes Q \simeq \mathcal{C}_{12}$,} so~$G$ would be abelian, a contradiction. Therefore,~\mbox{$[P,Q]\simeq \mathcal{C}_2\times \mathcal{C}_2$} and we conclude that~\mbox{$G\simeq \mathcal{C}_2 \times A_4$.}

\medskip

It remains to prove (iv). Clearly, $G/K \simeq A_4$. As a consequence, on one hand we have~\mbox{$Z(G)K/K \leq Z(G/K) = \{1\}$,} and so~$Z(G)\leq K$. On the other hand, since \mbox{$|G/K : C_G(K)/K| = |G:C_G(K)|\leq \lvert\operatorname{Aut}(K)\rvert = 2$,} we have~$C_G(K)=G$, and so \mbox{$K\leq Z(G)$.} Thus,~\mbox{$K=Z(G)$.} Let~$H/K$ denote the subgroup of order~$4$ in~$G/K$; in particular,~$H\trianglelefteq G$. Furthermore, let~$P$ be a Sylow $2$-subgroup of~$G$ contained in~$H$. Clearly, \mbox{$H\!=\!KP\!=\!Z(G)P$,} so $P$ is normal, and hence characteristic, in~$H$. Thus, \mbox{$P\trianglelefteq G$} and \mbox{$P\simeq H/K \simeq \mathcal{C}_2\times\mathcal{C}_2$.} Moreover, $G=P\rtimes Q$, where $Q\simeq\mathcal{C}_9$ is a Sylow $3$-subgroup of~$G$. Therefore,\linebreak \mbox{$G\simeq (\mathcal{C}_2\times\mathcal{C}_2)\rtimes \mathcal{C}_9$,} where the action is necessarily the unique non-trivial one.
\end{proof}

\begin{prop}\label{3177SylNormCycl}
Let $G$ be a finite group, and $p$ the largest prime dividing the order of~$G$. If $G$ is neither a $2$-group nor a $3$-group, and $\psi'(G) > \frac{31}{77}$, then $G$ has a normal cyclic~Sy\-low~\hbox{$2$-sub}\-group or a normal cyclic~Sy\-low~\hbox{$p$-sub}\-group.
\end{prop}
\begin{proof}
Let $|G|=n=p_1^{\alpha_1}p_2^{\alpha_2}\ldots p_k^{\alpha_k}$, where $k,\alpha_1,\alpha_2, \ldots, \alpha_k$ are positive integers and $p_1<p_2<\ldots<p_k=p$ are prime numbers. Throughout the proof, we will repeatedly use the following elementary observation: if $x\in G$ and $|G:\langle x \rangle| < p$, then~$\langle x \rangle$ contains a normal cyclic~Sy\-low~\hbox{$p$-sub}\-group of~$G$, because the order of $\mathrm{core}_G(\langle x \rangle)$ is coprime with~$p$.

Clearly, we may assume that $G$ is non-cyclic. Then, Lemma~\ref{smallestprime1943-3177}\ref{smallestprime3177} yields that the smallest prime dividing $|G|$ is either $2$ or $3$. Moreover, by hypothesis, it follows that~\hbox{$k\geq 2$.}

First, suppose $p_1=3$. Then, by Lemma~\ref{CyclBoundIndex}\ref{CyclBoundIndex_TWO}, there exists  $x\in G$ such that
$$|G:\langle x \rangle| < \frac{77}{31} \cdot \frac{4}{3} \cdot \frac{p+1}{p_2} < \frac{77}{31} \cdot \frac{4}{3} \cdot \frac{p+1}{5} < p \; ,$$
since $p_2\geq 5$ and the last inequality holds for every prime $p\geq 5$. Thus, $G$ has a normal cyclic~Sy\-low~\hbox{$p$-sub}\-group, and the statement is proved in this case.

Now, let us assume that $2$ divides the order of $G$. At this stage, the proof proceeds by a case-by-case analysis depending on the value of~$k$ and, if necessary, also on the value of~$p$. Let~\hbox{$K=\mathrm{core}_G(\langle x \rangle)$.}

\smallskip

    \begin{enumerate}[label={}, leftmargin=0pt]

    \item {\bf Case 1:} $k=2$.\\[3pt]
    By Lemma~\ref{CyclBoundIndex}\ref{CyclBoundIndex_ONE}, there exists $x\in G$ such that
    $$|G:\langle x \rangle| < \frac{77}{31}\cdot\frac{3}{2}\cdot\frac{4}{3} \approx 4.96\ldots \; .$$
   
    \begin{itemize}[leftmargin=15pt]
        \item Assume $p \geq 5$. Then $|G:\langle x \rangle| < p$, so $G$ has a normal cyclic Sylow $p$-subgroup.
        \item Assume $p=3$. Hence, $|G|=2^{\alpha_1} \cdot 3^{\alpha_2}$. We now examine the possible values of the index~\hbox{$|G:\langle x \rangle|$.}
        \smallskip
        \begin{itemize}[label=\footnotesize\ding{70}, leftmargin=15pt]
        \item Suppose $|G:\langle x \rangle| \leq 2$. Then $G$ has a normal cyclic Sylow $3$-subgroup.
        \smallskip
        \item Suppose $|G:\langle x \rangle| = 4$. By Lemma \ref{lucchini}, it follows that~\mbox{$|\langle x \rangle : K| < 4$.} If $\langle x \rangle \trianglelefteq G$ or~\hbox{$|\langle x \rangle:K|=2$,} then $K$ contains a normal cyclic Sylow $3$-subgroup of $G$. Therefore, we may assume that~\hbox{$|\langle x \rangle : K|=3$}, and in particular that~$\langle x \rangle$ is not normal in $G$. Consequently,~$|G:K|=12$; more precisely, $G/K\simeq A_4$. We now distinguish two scenarios, depending on whether 2 divides the order of $K$.
        \smallskip
        \begin{itemize}[label=\footnotesize\ding{71}, leftmargin=15pt]
        \item Assume that $2$ divides $|K|$. Then, there exists a subgroup $M$ of $K$ such that $|K:M|=2$. As a consequence, $G/M$ is a group of order $24$ with a non-normal subgroup of order $6$, namely $\langle x \rangle/M$, whose core in $G/M$, namely $K/M$, has order $2$.
        Hence, by Lemma \ref{techGAP}\ref{techGAP3}, we have $\psi(G/M)\leq 99$. Therefore, by a combination of Lemmas~\ref{BoundPsiCn}\ref{BoundPsiCn_TWO} and \ref{BoundPsiQuotient}, we obtain
        $$\frac{31}{77}\cdot\frac{2}{3+1}n^2 \leq \frac{31}{77}\psi(\mathcal{C}_n) < \psi(G) \leq \psi(G/M)\frac{n^2}{24^2} \leq \frac{99}{\, \, 24^2 \,}n^2 \, ,$$
         which yields a contradiction.
        \smallskip 
        \item Assume that $2$ does not divide $|K|$. In particular, $\alpha_1=2$. If $K=\{1\}$, we have $G\simeq G/K \simeq A_4$, so $\psi'(G)=\psi'(A_4)=\frac{31}{77}$, which contradicts the assumption. Thus, $3$ divides $K$. Then, there exists a subgroup $M$ of $K$ such that $|K:M|=3$. As a consequence, $G/M$ is a group of order $36$ with a non-normal subgroup of order $9$, namely $\langle x \rangle/M$, whose core in $G/M$, namely $K/M$, has order $3$.
        Hence, by Lemma \ref{techGAP}\ref{techGAP4}, we have $\psi(G/M)=265$. 
        Therefore, by a combination of Lemmas~\ref{psiCyclic} and \ref{BoundPsiQuotient}, we obtain
        $$\frac{31}{77}\cdot 11\cdot\frac{3^{2\alpha_2 +1}+1}{4} = \frac{31}{77}\psi(\mathcal{C}_n) < \psi(G) \leq \psi(G/M)|M|^2 = 265\cdot 3^{2\alpha_2 -4} \, ,$$
        which does not admit any real solutions, yielding a contradiction.
        \end{itemize}
        
\smallskip

        \item Suppose $|G:\langle x \rangle| = 3$. If $\langle x \rangle \trianglelefteq G$, then $G$ has a normal cyclic Sylow $2$-subgroup, so in the following we further assume that $\langle x\rangle$ is not normal in $G$. Thus, by Lemma~\ref{lucchini}, $|\langle x \rangle : K|=2$ and, consequently, \hbox{$|G:K|=6$.} We need to consider three distinct subcases, depending on the values of $\alpha_1$ and $\alpha_2$.
\smallskip
\begin{itemize}[label=\footnotesize\ding{71}, leftmargin=15pt]
    \item If $\alpha_1,\alpha_2\geq 2$, then there exists a subgroup $M$ of $K$ such that $|K:M|=6$. As a consequence, $G/M$ is a group of order $36$ with a non-normal cyclic subgroup of order $12$, namely $\langle x \rangle/M$. Hence, by Lemma \ref{techGAP}\ref{techGAP1}, we have $\psi(G/M)=243$. Therefore, by a combination of Lemmas~\ref{BoundPsiCn}\ref{BoundPsiCn_TWO} and \ref{BoundPsiQuotient}, we obtain
         $$\frac{31}{77}\cdot\frac{2}{3+1}n^2 \leq \frac{31}{77}\psi(\mathcal{C}_n) < \psi(G) \leq \psi(G/M)\frac{n^2}{36^2} = \frac{243}{\,\, 36^2 \,}n^2 \, ,$$
         which yields a contradiction.

         \smallskip

    \item If $\alpha_1=1$, then $|G|=2\cdot 3^{\alpha_2}$. For $\alpha_2=1$, we have $G\simeq\mathbb{Z}_6$ or $G\simeq S_3$, and in both cases $G$ has a normal cyclic Sylow $3$-subgroup. Suppose $\alpha_2 = 2$, so that~$|G|=18$. Therefore, $G=P\rtimes Q$, where $P$ is a Sylow $3$-subgroup of $G$ and $Q\simeq\mathcal{C}_2$ is a~\hbox{Sy}\-low $2$-subgroup of $G$. If $P\simeq\mathcal{C}_9$, then $G$ has a normal cyclic $3$-subgroup, so we instead assume that $P\simeq\mathcal{C}_3\times\mathcal{C}_3$. Let's compute $\psi(G)$. The identity has order $1$. In $P\setminus\{1\}$ we find $8$ elements of order $3$. Since $\langle x \rangle$ is not normal in~$G$, we have that its conjugacy class has size $3$, so there are at least $3$ subgroups of order $6$ and, consequently, at least $6$ elements of order $6$. Similarly, since $Q$ is not normal in $G$, we have at least $3$ elements of order $2$. We have accounted for~$18$ elements, so $\psi(G) = 1 + 8\cdot 3 + 6\cdot 6 + 3\cdot 2 = 67$. Hence, $\psi'(G) = \frac{67}{183} < \frac{31}{77}$, which is a contradiction.
        
    \quad Now, assume furthermore that $\alpha_2\geq 3$. Then, there exists a subgroup $M$ of~$K$ such that $|K:M|=9$. As a consequence, $G/M$ is a group of order $54$ with a non-normal cyclic subgroup of order $18$, namely $\langle x \rangle/M$. Hence, by~Lem\-ma~\ref{techGAP}\ref{techGAP2}, we have $\psi(G/M)=553$. 
        
        Therefore, by a combination of Lemmas~\ref{BoundPsiCn}\ref{BoundPsiCn_TWO} and~\ref{BoundPsiQuotient}, we obtain
         $$\frac{31}{77}\cdot\frac{2}{3+1}n^2 \leq \frac{31}{77}\psi(\mathcal{C}_n) < \psi(G) \leq \psi(G/M)\frac{n^2}{54^2} = \frac{553}{\, \, 54^2 \,}n^2 \, ,$$
         which yields a contradiction.
         
         \smallskip
         
        \item If $\alpha_2=1$, then $|G|=2^{\alpha_1}\cdot 3$. Let $P$ be a Sylow \mbox{$3$-subgroup} of~$G$. Clearly,~\mbox{$P\simeq\mathcal{C}_3$.} Since $|G:K|=6$, it follows that $|K|=2^{\alpha_1-1}$. Thus,~\mbox{$KP=K\rtimes P= K\times P$} is abelian. Moreover, since \mbox{$|G:KP|=2$},  \mbox{$KP\trianglelefteq G$.} Therefore, $P\trianglelefteq G$. Hence, $G$ has a normal cyclic Sylow $3$-subgroup.
\end{itemize}
         
    \end{itemize}
    \end{itemize}

\smallskip
    
        \item {\bf Case 2:} $k = 3$.\\[3pt]
        By Lemma \ref{CyclBoundIndex}\ref{CyclBoundIndex_ONE}, there exists $x\in G$ such that
        $$|G:\langle x \rangle| < \frac{77}{31} \cdot \frac{3}{2} \cdot \frac{4}{3} \cdot \frac{6}{5} \approx 5.96\ldots \, .$$
     
        \begin{itemize}[leftmargin=15pt]
        \item Assume $p\geq 7.$ Then $|G:\langle x \rangle| < p$, and so $G$ has a normal cyclic Sylow $p$-subgroup.
\smallskip
        \item Assume $p=5$. First, if $|G:\langle x \rangle|<5=p$, then $G$ has a normal cyclic Sylow \mbox{$p$-subgroup.} Next, suppose that $|G:\langle x \rangle|=5$.  By Lemma \ref{lucchini}, \mbox{$|\langle x \rangle:K| < |G:\langle x \rangle| = 5$.} If~\mbox{$|\langle x \rangle:K|\in\{1,3\}$}, then $G$ has a normal cyclic Sylow $2$-subgroup. Suppose this is not the case, so that $|\langle x \rangle:K|\in\{2,4\}$. Then $G$ has a normal cy\-clic~Sy\-low~\hbox{$3$-sub}\-group, say~$P_3$. Thus, there exists a subgroup~$H$ such that $G=P_3\rtimes H$; in particular~\mbox{$|G|=|P_3||H|$.} Then, clearly, we can write $H=\langle P_2, P_5\rangle$, where $P_2$ and $P_5$ are, respectively, a
        Sylow $2$-subgroup  and a~Sy\-low~\hbox{$5$-sub}\-group of $G$. Since~$\pi(\operatorname{Aut}(P_3))\subseteq\{2,3\}$, it follows that $[P_3,P_5] = \{1\}$, hence $P_3\leq N_G(P_5)$. 
        By~Lem\-ma~\ref{BoundPsiPNormalCyclicSylow}\ref{BoundPsiPNormalCyclicSylowBound}, we have
        $$\frac{31}{77} < \psi'(G) \leq \psi'(P_3)\psi'(G/P_3)=\psi'(H) \; .$$ 
        It is therefore legimate to apply the case $k=2$ to the group $H$, thereby obtaining that~$H$ has a cyclic normal Sylow $s$-subgroup, where $s=2$ or $s=5$. If~$s=5$, then~\mbox{$P_2 \leq N_G(P_5)$,} and so $P_5$ is a normal cyclic Sylow subgroup of $G$. Now, consider the case $s=2$. Since $\operatorname{Aut}(P_2)$ is a $2$-group, it follows that $[P_2,P_5]=\{1\}$, and hence~\mbox{$P_5\leq Z(G)$;} in particular, $G = P_5 \times (P_3 \rtimes P_2)$. Therefore, if we assume that $P_5$ is not cyclic, then the successive application of the multiplicativity of $\psi'$, Lem\-ma~\ref{BoundPsiPNormalCyclicSylow}\ref{BoundPsiPNormalCyclicSylowBound} and~The\-o\-rem~\ref{BoundPsiNonCycl} leads to a contradiction:
        $$ \frac{31}{77} < \psi'(G) = \psi'(P_5)\psi'(P_3\rtimes P_2) \leq  \psi'(P_5)\psi'(P_3)\psi'(P_2) = \psi'(P_5) < \frac{1}{4} \; .$$
        Thus, once again, $P_5$ turns out to be a normal cyclic Sylow subgroup of $G$. Therefore, for both choices of $s$, $G$ has a normal cyclic Sylow $p$-subgroup.
        
        \hspace{0.4cm} Note that by using Burnside's normal $p$-complement theorem one immediately obtains normality of $P_5$ in $G$, slightly simplifying some arguments, but we chose not to do so in order to keep the argument more self-contained.
        
        \end{itemize}
        
\smallskip
        
        \item {\bf Case 3:} $k\geq 4$.\\[3pt]
        By Lemma \ref{CyclBoundIndex}\ref{CyclBoundIndex_TWO}
, there exists $x\in G$ such that
    $$|G:\langle x \rangle| < \frac{77}{31} \cdot \frac{3}{2} \cdot \frac{4}{3} \cdot\frac{6}{5} \cdot \frac{p+1}{7} < p \; ,$$
    since the last inequality is true for every prime $p \geq 7$. Hence, $G$ has a normal cyclic~Sy\-low~\hbox{$p$-sub}\-group.
    \end{enumerate}
\smallskip

The proof is complete. 
\end{proof}

\begin{coroll}\label{1943SylNormCycl}
Let $G$ be a finite group, and $p$ the largest prime dividing the order of~$G$. If $G$ is not a $2$-group, and $\psi'(G) \geq \frac{19}{43}$,  then $G$ has a normal cyclic Sylow $2$-subgroup or a normal cyclic Sylow $p$-subgroup.
\end{coroll}
\begin{proof} $G$ is not a $2$-group by hypothesis, nor is it a $3$-group by Lemma~\ref{smallestprime1943-3177}\ref{smallestprime1943}. Moreover, $\psi'(G) \geq \frac{19}{43} > \frac{31}{77}$. Therefore, by Proposition~\ref{3177SylNormCycl}, the claim follows.
\end{proof}

As a further step towards the proof of Theorem A, we focus on the special case of \hbox{$M$-groups.} To this end, we first obtain an explicit formula for the $\psi'$-value of a restricted class of $P^*$-groups, and subsequently examine its behavior within the range of interest.

\begin{prop}\label{PsiP*}
Let $p$ be an odd prime number and $n,k$ positive integers. Then
$$\psi'(P^*(p^n,2^k)) = \frac{p+1}{p^{2n+1}+1} + \frac{(p+1)(p^n-1)}{p^{2n+1}+1}\cdot \frac{(p+3)2^{2k-1}+p}{2^{2k+1}+1}  \ .$$
\end{prop}
\begin{proof}
    Let $G=P^*(p^n, 2^k)$, so $G = A \rtimes_\iota \langle x\rangle$, where $A$ is an elementary abelian group of order $p^n$ and $|\langle x\rangle| = 2^k$; in particular, every element of $G$ can be written \textit{uniquely} as $ax^m$, where $a\in A$ and $m$ is an integer such that $0\leq m \leq 2^k - 1$ (this uniqueness establishes a natural bijection $ax^m\in G \mapsto (a,m)\in A\times\{0,1,\dots,2^k-1\}$,  which underlies our counting). By means of the parameters $a$ and $m$, we partition $G$ into the following four subsets:\\
$$\begin{array}{l}
S_1\coloneqq \big\{ \, ax^m \, : \, a=1, \, m\in\{0,1,\ldots, 2^k-1\} \, \big\} = \langle x \rangle \; ,\\[8pt]
S_2 \coloneqq \big\{ \, ax^m \, : \, a\in A\setminus\{1\}, \, m=0 \, \big\} = A\setminus\{1\} \; ,\\[8pt]
S_3 \coloneqq \big\{ \, ax^m \, : \, a\in A\setminus\{1\}, \, m\in\{1,3,\ldots, 2^k-1\} \, \big\} \; ,\\[8pt]
S_4 \coloneqq \big\{ \, ax^m \, : \, a\in A\setminus\{1\}, \, m\in\{2,4,\ldots, 2^k-2\} \, \big\} \; .
\end{array}
$$
\ \\
For each of them, we compute the corresponding value of $\psi$, that is, the sum of the orders of the elements of the subset.
\begin{itemize}
[leftmargin=15pt, listparindent=\parindent, parsep=0pt]
\item Clearly, by Lemma \ref{psiCyclic}, $$\psi(S_1) = \frac{2^{2k+1}+1}{3} \; .$$
\item Since $S_2=A\setminus\{1\}$ has $p^n-1$ elements, each of order $p$, we have
$$\psi(S_2) = p(p^n-1) \; .$$
\item
\begin{sloppypar}
 The elements of $S_3$ are, by definition, those of the form $ax^m$, with~\mbox{$a\in A\setminus\{1\}$} and \hbox{$m\in\{1,3,\ldots,2^k-1\}$.} Since there are \hbox{$p^n-1$} choices for $a$, and \hbox{$|\{1,3,\ldots,2^k -1 \}| = 2^{k-1}$} choices for $m$, we have exactly~$(p^n-1)\cdot 2^{k-1}$ elements in $S_3$. Moreover, since $m$ is odd, a direct computation shows that all these elements share the same order $2^k$. Therefore,
$$\psi(S_3) = (p^n-1)\cdot 2^{k-1} \cdot 2^k = 2^{2k-1}(p^n-1) \; .$$
\end{sloppypar}
\item We now focus on the elements $ax^m\in G$ with $a\in A\setminus\{1\}$ and $m\in\{2,4,\ldots,2^k-2\}$, which constitute the set $S_4$. Assume $k\geq 2$, so that $S_4$ is non-empty. Each $m$, being even, can be written \textit{uniquely} as $m = 2^s d$, where $s$ is a positive integer and $d$ is an odd positive integer; more precisely, since $m\in\{2,4,\dots, 2^k - 2\}$, we have $1\leq s \leq k-1$ and~$d\in \{1,3,\ldots,2^{k-s}-1\}$. Indeed, the map $ax^m=ax^{2^s d} \mapsto (a,s,d)$ describes a natural bijection between the elements of $S_4$ and all triples $(a,s,d)$ satisfying the respective constraints on $a,s,d$.

\begin{sloppypar}
Now, let $s$ be fixed and restrict our attention to the elements of the form~$ax^{2^s d}$. There are~\hbox{$p^n - 1$} choices for $a$, and $|\{1,3,\ldots,2^{k-s}-1\}|=2^{k-s-1}$ choices for~$d$, giving a total of~\hbox{$(p^n - 1)\cdot 2^{k-s-1}$} such elements. Moreover, since~$[a,x^2]=1$, we have~\hbox{$o(ax^{2^s d}) = \operatorname{lcm}\!\big(o(a),o(x^{2^s d})\big) =  2^{k-s}p$} ($p$ is odd by assumption). Summarizing, for each fixed $s$, there are exactly~$(p^n - 1)\cdot 2^{k-s-1}$ elements of the form $ax^{2^s d}$, all sharing the same order $2^{k-s}p$. Therefore, summing over all possible values of $s$, we obtain~$\psi(S_4)$ as follows: 
\end{sloppypar}
\vspace{-0.5cm}
$$\begin{array}{ccl}
\displaystyle \psi(S_4) &\displaystyle = & \displaystyle \sum_{s=1}^{k-1} \psi\!\ \Big(\big\{\, ax^{2^s d} \, : \, a\in A\setminus\{1\}, \ d\in\{1, 3, \ldots, 2^{k-s}-1 \} \, \big\}\Big) \ = \\[18pt]
\ &\displaystyle = & \displaystyle \sum_{s=1}^{k-1} (p^n - 1)\cdot 2^{k-s-1} \cdot 2^{k-s}p \ = \ p(p^n-1) 2^{2k-1} \sum_{s=1}^{k-1} (2^{-2})^s \ = \\[18pt]
\ &\displaystyle = & \displaystyle p(p^n-1) 2^{2k-1} \left[ \frac{(2^{-2})^{k}-1}{2^{-2}-1} - 1 \right] \ = \ p(p^n-1) \frac{2^{2k-1}-2}{3} \; .
\end{array}$$
Note that this expression also covers the case $k=1$, for which $S_4$ is empty, consistently yielding $0$.
\end{itemize}

\noindent Clearly, $\psi(G)=\psi(S_1)+\psi(S_2)+\psi(S_3)+\psi(S_4)$. Finally, dividing $\psi(G)$ by 
    $$\psi(\mathcal{C}_{p^n2^k}) = \frac{p^{2n+1}+1}{p+1}\cdot \frac{2^{2k+1}+1}{3}$$ (see Lemma \ref{psiCyclic}), we arrive at the stated expression for $\psi'(G)$.
\end{proof}

\begin{coroll}\label{P*(2)1943}
    Let $p$ be an odd prime number, and let $n,k$ be positive integers. Then\textnormal:
    $$\psi'\big(P^*(p^n,2^k)\big) > \frac{19}{43} \mbox{ if and only if either } p=3, n=1,k\geq1 \mbox{ or } p=5, n=1, k\leq 2.$$
\end{coroll}

\begin{proof}
    By substituting $p=3$ and $n=1$ into the expression given by Lemma \ref{PsiP*}, it is easy to verify that 
    $$\psi'(P^*(3,2^k)) = \frac{2^{2k+3} + 7}{7 \cdot 2^{2k+1}+7} > \frac{19}{43} \; $$ for all $k\geq 1$.
    Similarly, by substituting $p=5$, $n=1$, and either $k=1$ or $k=2$, we obtain, respectively, 
    $$\psi'(P^*(5,2)) = \frac{31}{63} > \frac{19}{43}\mbox{\quad and \;\:} \psi'(P^*(5,4)) = \frac{103}{231} > \frac{19}{43} \; .$$ 
 
Let us now consider the more difficult implication. We begin by pointing out the following monotonicity properties that will be used in the ensuing argument. First, the functions
    $$\frac{(p+3)2^{2k-1}+p}{2^{2k+1}+1} \; , \quad \frac{p}{p^{2n+1}+1} \; , \quad \frac{(p+1)(p^n-1)}{p^{2n+1}+1}$$
are decreasing in $k$, $n$ and $n$, respectively. Second, all single-variable functions (namely,~$p$ or $k$) appearing in what follows, as an entire side of an inequality, are decreasing. 
    
Now, if $n\geq 2$, we obtain
    $$\begin{array}{ccl}
    \displaystyle \psi'(P^*(p^n,2^k)) & \overset{\ref{PsiP*}}{=} & \displaystyle \frac{p+1}{p^{2n+1}+1} + \frac{(p+1)(p^n-1)}{p^{2n+1}+1}\cdot \frac{(p+3)2^{2k-1}+p}{2^{2k+1}+1} \\[20pt]
    \displaystyle & \displaystyle \overset{k=1}{\leq} & \displaystyle \frac{p+1}{p^{2n+1}+1} + \frac{(p+1)(p^n-1)}{(p^{2n+1}+1)}\cdot\frac{p+2}{3} \\[20pt]
    \displaystyle & \overset{n=2}{\leq} & \displaystyle \frac{(p+1)(p^3+2p^2-p+1)}{3(p^5+1)} \; \overset{p=3}{\leq} \; \frac{43}{183} < \frac{19}{43} \; ,
    \end{array}$$
which is a contradiction. Hence, $n=1$ and
$$ \displaystyle \psi'(P^*(p^n,2^k)) = \displaystyle \frac{p+1}{p^3+1} + \frac{p^2-1}{p^3+1}\cdot \frac{(p+3)2^{2k-1}+p}{2^{2k+1}+1} \; .$$

Next, if $p\geq 7$, then
    $$ \displaystyle \psi'(P^*(p^n,2^k)) \overset{k=1}{\leq} \displaystyle \frac{p^2+p+1}{3p^2-3p+3} \; \overset{p=7}{\leq} \; \frac{19}{43} \; ,$$
resulting again in a contradiction, so $p=3$ or $p=5$.

\smallskip
    
Finally, in the case $p=5$ we have to add the restriction $k\leq 2$, otherwise 
    $$\displaystyle \psi'(P^*(p^n,2^k)) =
    \displaystyle \frac{3\cdot 2^{2k+1}+7}{7\cdot 2^{2k+1}+7} \displaystyle \; \overset{k=3}{\leq} \; \displaystyle \; \frac{391}{903} < \frac{19}{43} \; ,$$
which once more contradicts the assumption.

\smallskip
    
We can now conclude that the values of $p$, $n$, $k$ are restricted as stated.
\end{proof}

The following proposition shows that, by exploiting the well-known structure of finite $M$-groups, the classification described in Theorem A can be obtained quite easily. Later, the assumption of being an $M$-group will turn out to be dispensable.

\begin{prop}\label{Mg1943}
    Let $G$ be a finite $M$-group with $\psi'(G) > \frac{19}{43} $. Then, up to isomorphism,~$G$ belongs to one of the following pairwise disjoint classes of groups\textnormal:
    \begin{itemize}
        \item cyclic groups;
        \item $\left(\mathcal{C}_{2^{n-1}}\!\times \mathcal{C}_2\right) \times \mathcal{C}_m$, where $n\geq 2$ and $(m,2)=1$;
        \item $M(2^n) \times \mathcal{C}_m$, where $n\geq 4$ and $(m,2) =1$;
        \item $Q_8 \times \mathcal{C}_m$, where $(m,2)=1$;
        \item $\left(\mathcal{C}_3 \rtimes_{\iota} \mathcal{C}_{2^k}\right) \times \mathcal{C}_m$, where $k\geq 1$ and $(m,6)=1$;
        \item $\left(\mathcal{C}_5 \rtimes_{\iota} \mathcal{C}_2\right) \times \mathcal{C}_m$, where $(m,10)=1$;
        \item $\left(\mathcal{C}_5 \rtimes_{\iota} \mathcal{C}_4\right) \times \mathcal{C}_m$, where $(m,10)=1$.
    \end{itemize}
\end{prop}

\begin{proof}
Since $G$ is an $M$-group, by Theorem \ref{StructureFiniteMgroups} we have the decomposition
$$G = P_1 \times \ldots \times P_t \; ,$$
where the factors have pairwise coprime orders, and each of them is either a prima\-ry~\hbox{$M$-group} or a $P^*$-group. Suppose, for contradiction, that there are at least two non-cyclic~$P_i$'s. Then, by Theorem \ref{BoundPsiNonCycl}, we have
$$\frac{19}{43} < \psi'(G) = \psi'(P_1)\cdot\ldots\cdot\psi'(P_t) \leq \left( \frac{7}{11} \right)^2 \; ,$$
\begin{sloppypar}
\noindent a contradiction. Hence, there is at most one non-cyclic $P_i$, say $P_1$, and consequent\-ly~\hbox{$\psi'(P_1)>\frac{19}{43}$.} If $P_1$ is cyclic, then $G$ is cyclic. Suppose that $P_1$ is non-cyclic. By~Lem\-ma~\ref{smallestprime1943-3177}\ref{smallestprime1943}, $2$ divides the order of $P_1$, so $P_1$ is either a non-cyclic $2$-primary $M$-group or a $P^*(p^n, 2^k)$-group for some odd prime $p$ and certain positive integers $n,k$.  Applying~Proposition \ref{2g-classification}\ref{2g>1943} in the former case, and Corollary \ref{P*(2)1943} in the latter case, we conclude that, up to isomorphism, $G$ belongs to one of the classes of groups listed in the statement.
\end{sloppypar}
\end{proof}

We are now in a position to establish Theorem A. However, before proceeding with the proof, we present a few elementary technical lemmas, which will serve to streamline certain recurring computations.

\begin{lemma}\label{BoundP/psi(C_P)}
    Let $p$ be a prime number, $n$ a positive integer and $P$ a group of order $p^n$. The following hold\textnormal:
    \begin{enumerate}[label=$(\roman*)$]
        \item If $p\geq 3$, then $\displaystyle \frac{|P|}{\psi(\mathcal{C}_{|P|})} \leq \frac{3}{7} \;$.
        \label{B1}
        \item If $p\geq 5$, then $\displaystyle \frac{|P|}{\psi(\mathcal{C}_{|P|})} \leq \frac{5}{21} \;$.
        \label{B2}
        \item If $p\geq 7$, then $\displaystyle \frac{|P|}{\psi(\mathcal{C}_{|P|})} \leq \frac{7}{43} \;$.
        \label{B3}
        \item If $p\geq 11$, then $\displaystyle \frac{|P|}{\psi(\mathcal{C}_{|P|})} \leq \frac{11}{111} \;$.
        \label{Bplus}
        \item If $n\geq 2$ and $p\geq3$, then $\displaystyle \frac{|P|}{\psi(\mathcal{C}_{|P|})} \leq \frac{9}{61} \;$.
        \label{B4}
    \end{enumerate}
\end{lemma}
\begin{proof}
    By Lemma \ref{psiCyclic}, it is immediate that
    $$\frac{|P|}{\psi(\mathcal{C}_{|P|})} = \frac{p^n(p+1)}{p^{2n+1}+1} \; .$$
    Since the latter function is monotonically decreasing in both variables over the domain of interest ($n\geq 1$, $p\geq 2$), we obtain the listed bounds by evaluating it at the points $(n,p) = (1,3), (1,5), (1,7), (1,11), (2,3)$.
\end{proof}

\begin{lemma}\label{techCycl}
    Let $m$ be a positive integer strictly greater than $1$, and let $p$ and $q$ be prime numbers dividing $m$. If $q<p$, then $\psi(\mathcal{C}_{m/q}) > \psi(\mathcal{C}_{m/p})$. In particular, if $q$ is the smallest prime dividing $m$, then $\mathcal{C}_{m/q}$ maximizes the value of $\psi$ among the proper subgroups of $\mathcal{C}_m$.
\end{lemma}
\begin{proof}
Let $q^\alpha$ and $p^\beta$ be, respectively, the highest powers of $q$ and $p$ that divide $m$.
By~Lemma~\ref{psiCyclic}, we obtain that our inequality is equivalent to
$$\frac{q^{2\alpha-1}+1}{q+1}\cdot\frac{p^{2\beta+1}+1}{p+1} > \frac{q^{2\alpha+1}+1}{q+1}\cdot\frac{p^{2\beta-1}+1}{p+1} \; ,$$
that is,
$$q^{2\alpha-1}(p^{2\beta+1}+1)+p^{2\beta+1} > q^{2\alpha-1}q^2(p^{2\beta-1}+1)+p^{2\beta-1} \; .$$
Clearly, $p^{2\beta+1}>p^{2\beta-1}$. Moreover, since $p\geq q+1$, we have $$p^{2\beta+1}+1 > p^2p^{2\beta-1} \geq (q+1)^2p^{2\beta-1} = q^2p^{2\beta-1}+(2q+1)p^{2\beta-1} \geq q^2(p^{2\beta-1}+1).$$ Thus, the inequality is proved. The last part of the statement follows immediate\-ly~from the fact that every maximal subgroup of $\mathcal{C}_m$ is isomorphic to $\mathcal{C}_{m/r}$ for some prime~$r$.
\end{proof}

\begin{lemma}\label{tech1}
Consider the direct product $A\times \mathcal{C}_m$, where $A$ is a finite group and $m$ is a positive odd integer such that $(|A|,m)=1$. Let $B$ be a maximal subgroup of $A$ that attains the maximum value of $\psi$ among the maximal subgroups of $A$. The following statements hold.
    \begin{enumerate}[label=$(\roman*)$]
        \item If \kern0.12em \scalebox{0.8}{$\dfrac{\psi(A)}{\psi(B)}$}\kern0.25em$\leq 7$, then $B\times \mathcal{C}_m$ is a maximal subgroup of $A\times\mathcal{C}_m$ that attains the maximum value of~$\psi$ among the maximal subgroups \textnormal(and hence among all proper subgroups\textnormal) of~\hbox{$A\times\mathcal{C}_m$.}
        \label{tech1.1}
        \item If $B$ is cyclic and $|A:B|=2$, then \kern0.1em \scalebox{0.8}{$\dfrac{\psi(A)}{\psi(B)}$}\kern0.25em$\leq 7$, so the conclusion of point $(i)$ follows.
        \label{tech1.2}
    \end{enumerate}
\end{lemma}
\begin{proof}
$(i)$ \, For $m=1$ there is nothing to prove, so we assume $m>1$. Since the orders of~$A$ and $\mathcal{C}_m$ are coprime, 
    the maximal subgroups of $A\times \mathcal{C}_m$ are precisely those of the form $H\times K$, where either $H=A$ and $K$ is a maximal subgroup of $\mathcal{C}_m$, or $H$ is a maximal subgroup of $A$ and $K=\mathcal{C}_m$. As $\psi(H\times K)=\psi(H)\psi(K)$, the maximum value of $\psi$ among the proper subgroups of $A\times \mathcal{C}_m$ is attained at either $B\times \mathcal{C}_m$ or $A\times \mathcal{C}_{m/q}$, where $q$ is the smallest prime dividing $m$ (see Lemma \ref{techCycl}). Let $\alpha$ be the exponent of the highest power of $q$ dividing $m$. Then, by hypothesis, 
    $$\frac{\psi(A)}{\psi(B)} \leq 7 \leq \frac{q^{2\alpha+1}+1}{q^{2\alpha-1}+1} = \frac{\psi(\mathcal{C}_m)}{\psi(\mathcal{C}_{m/q})} \; .$$
    Note that this latter bound holds because the right-hand side function is increasing with respect to both variables over the domain of interest ($q\geq 3$, $\alpha\geq 1$), and so attains its global minimum value, which is $7$, at the point $(q,\alpha)=(3,1)$. In conclusion, it follows that $\psi(A\times \mathcal{C}_{m/q}) \leq \psi(B\times \mathcal{C}_m)$.

\medskip

\noindent$(ii)$ \, Let $\alpha$ be the exponent of the highest power of $2$ dividing $|A|$. By the hypotheses, we~have
    $$\frac{\psi(A)}{\psi(B)} \leq \frac{\psi(\mathcal{C}_{|A|})}{\psi(\mathcal{C}_{|B|})} = \frac{2^{2\alpha+1}+1}{2^{2\alpha -1}+1} < \lim_{\alpha\to+\infty}\frac{2^{2\alpha+1}+1}{2^{2\alpha -1}+1}=4 < 7 \, .$$ 

\ \vspace{0.05cm}

\noindent The statement is proved.
\end{proof}

\bigskip

\noindent {\bf Proof of Theorem A.}\\
We proceed by induction on the order of $G$ to prove that $G$ is an $M$-group.
Clearly, $\{1\}$ is an $M$-group. Assume $|G|>1$, and that every group of smaller order than $G$ with $\psi'$-value exceeding $\frac{19}{43}$ is an $M$-group.

Clearly, we may assume that $G$ is not cyclic, so in particular $2$ divides the order of~$G$ by~Lem\-ma~\ref{smallestprime1943-3177}\ref{smallestprime1943}, and moreover, by Proposition \ref{2g-classification}\ref{2g>1943}, we may assume that $G$ is not a $2$-group.
By~Corollary~\ref{1943SylNormCycl}, there exists a normal cyclic Sylow $r$-subgroup $P$ of $G$, for some \mbox{$r\in\{2,p\}$,} where~$p$ is the largest prime dividing the order of $G$ (evidently, $p\geq 3$).
Thus, by the~\hbox{Schur--Zas}\-sen\-haus theorem,
$$G = P \rtimes H \; ,$$
where $H$ is a suitable subgroup of $G$; in particular, $(|P|,|H|)=1$.

 Applying Lemma \ref{BoundPsiPNormalCyclicSylow}\ref{BoundPsiPNormalCyclicSylowBound}, we easily have that
$$\frac{19}{43} < \psi'(G)\leq \psi'(P)\psi'(G/P) = \psi'(G/P) \, .$$
Thus, by induction hypothesis, we obtain that $H \simeq G/P$ is an $M$-group with $\psi'(H)>\frac{19}{43}$, so $H$ must necessarily be one of the groups listed in the statement of Proposition~\ref{Mg1943}. After establishing some necessary preliminaries, we will analyze all the possibilities for $H$.

First, note that if $[P,H]=\{1\}$, which is equivalent to saying that $G$ is not merely a semidirect product of $P$ by $H$ but actually their direct product, that is, 
$$G = P\times H \; ,$$
then $G$ is an $M$-group, and the claim follows. Thus, we may assume that $[P,H]\neq\{1\}$. Consequently, $r=p$ (otherwise $r=2$, so $\operatorname{Aut}(P)$ would be a $2$-group, and we would again have $[P,H]=\{1\}$); in particular, $2$ divides the order of~$H$.

Next, we will need a suitable upper bound for $\psi'(G)$. From this point forward, given a finite group $A$, we denote by $l(A)$ the maximum value of $\psi$ among the maximal subgroups of $A$ (equivalently, among all the proper subgroups of $A$). Since $[P,H]\neq\{1\}$, we have \hbox{$C_H(P) < H$,} and hence
$$ \frac{\psi(C_H(P))}{\psi(\mathcal{C}_{|H|})} \leq \frac{l(H)}{\psi(\mathcal{C}_{|H|})} \; .$$
Starting from the equality in Lemma~\ref{BoundPsiPNormalCyclicSylow}\ref{BoundPsiPNormalCyclicSylowFormulaSDP}, dividing both sides by $\psi(\mathcal{C}_{|G|})$, and using the inequality above, we obtain
\begin{equation}\label{stareq}
\psi'(G) \leq \frac{|P|}{\psi(\mathcal{C}_{|P|})} \cdot \psi'(H) + \left(1-\frac{|P|}{\psi(\mathcal{C}_{|P|})} \right) \cdot \frac{l(H)}{\psi(\mathcal{C}_{|H|})} \ . \tag{$\star$} 
\end{equation}
Inequality $\eqref{stareq}$ is the bound we were aiming for.

We now proceed with the analysis of all possible cases for $H$, each time evaluating~$\psi'(H)$ and $l(H)$ --- with the aid of Lemma~\ref{techCycl} or Lemma~\ref{tech1} --- in order to exploit the bound provided by inequality~$\eqref{stareq}$, and subsequently making use of Lemma~\ref{BoundP/psi(C_P)} along with certain monotonicity properties of the involved functions. In each case, we will either conclude that $G$ is an $M$-group or reach the contradiction $\psi'(G) \leq \frac{19}{43}$.

\bigskip

\noindent {\bf Case 1:} {\it $H\simeq \mathcal{C}_m$, for some positive integer $m$}.

\smallskip

Denote by $n$ the positive integer such that $|P|=p^n$, and let $k$ be the exponent of the highest power of 2 dividing $m$. Since $2$ divides $|H|=m$, we have $k\geq 1$, and Lemma~\ref{techCycl} yields $l(H)=\psi(\mathcal{C}_{m/2})$. Moreover, clearly, $\psi'(H)=1$.
Therefore,
$$\begin{array}{ccl}
    \displaystyle \psi'(G) & \displaystyle \overset{\eqref{stareq}}{\leq} & \displaystyle\frac{|P|}{\psi(\mathcal{C}_{|P|})}\cdot 1 + \left(1-\frac{|P|}{\psi(\mathcal{C}_{|P|})} \right)\cdot\frac{2^{2k-1}+1}{2^{2k+1}+1} \;  \\[20pt]
    \displaystyle & \displaystyle  \overset{k=1}{\leq}  & \displaystyle\frac{|P|} {\psi(\mathcal{C}_{|P|})} \; + \left(1-\frac{|P|}{\psi(\mathcal{C}_{|P|})} \right)\cdot\frac{1}{3} = \frac{|P|} {\psi(\mathcal{C}_{|P|})}\cdot\frac{2}{3} + \frac{1}{3} \; .
    \end{array}$$
    Hence, if $p\geq 7$, then
    $$\displaystyle \psi'(G) \; \leq \; \frac{|P|}{\psi(\mathcal{C}_{|P|})}\cdot\frac{2}{3} + \frac{1}{3} \; \overset{\text{\ref{BoundP/psi(C_P)}\ref{B3}}}{\leq} \; \displaystyle\frac{7}{43}\cdot \frac{2}{3} + \frac{1}{3} = \frac{19}{43} \; ,$$
    and if $n\geq 2$, then
    $$\displaystyle \psi'(G) \; \leq \; \frac{|P|}{\psi(\mathcal{C}_{|P|})}\cdot\frac{2}{3} + \frac{1}{3} \; \overset{\text{\ref{BoundP/psi(C_P)}\ref{B4}}}{\leq} \; \frac{9}{61}\cdot\frac{2}{3} + \frac{1}{3}= \frac{79}{183} < \frac{19}{43} \; ,$$
    so, in both cases we reach a contradiction. Thus, $p=3$ or $p=5$, and $n=1$. As a result, since $p$ is the maximum of $\pi(G)$, we obtain that $|P|=3$ and $H$ is a $2$-group, or $|P|=5$ and $\pi(H)\subseteq\{2,3\}$.
    
    If $|P|=3$, then the only non-trivial way in which $H$ can act on~$P\simeq\mathcal{C}_3$ is by inversion. Hence, $G$ is isomorphic to $\mathcal{C}_3 \rtimes_{\iota} \mathcal{C}_{2^k}\simeq P^*(3,2^k)$, which is an $M$-group.
    
    Next, suppose $|P|=5$. Clearly, $H = H_2 \times H_3$, where $H_2$ and $H_3$ denote, respectively, the Sylow 2-subgroup and the Sylow 3-subgroup of $H$; in particular, $G=P\rtimes\left(H_2\times H_3 \right)$. Since $\operatorname{Aut}(P)\simeq\mathcal{C}_4$, $H_3$ acts trivially on $P$, so $[H_3,P]=\{1\}$, and hence $H_3\leq Z(G)$. Thus, $G=H_3\times \left( P \rtimes H_2 \right)$. Therefore,
    $$\frac{19}{43} < \psi'(G) = \psi'(H_3) \psi'(PH_2) = \psi'(PH_2) \; .$$
    Consequently, if $H_3\neq\{1\}$, then $PH_2$ is an $M$-group by the induction hypothesis, and so is $G$; in particular, we may assume $H_3=\{1\}$. Moreover, if $k\geq 3$, then
    $$\begin{array}{ccl}
    \displaystyle \psi'(G) & \displaystyle\overset{\eqref{stareq}}{\leq} & \displaystyle\frac{|P|}{\psi(\mathcal{C}_{|P|})} \cdot 1 + \left(1-\frac{|P|}{\psi(\mathcal{C}_{|P|})} \right)\cdot\frac{2^{2k-1}+1}{2^{2k+1}+1}  \\[20pt]
    \displaystyle & \displaystyle \overset{k=3}{\leq} & \displaystyle\frac{|P|}{\psi(\mathcal{C}_{|P|})} + \left(1-\frac{|P|}{\psi(\mathcal{C}_{|P|})} \right)\cdot\frac{11}{43} \; =
    \displaystyle\frac{|P|} {\psi(\mathcal{C}_{|P|})}\cdot\frac{32}{43} + \frac{11}{43}  \\[20pt] 
    \displaystyle & \displaystyle \overset{\text{\ref{BoundP/psi(C_P)}\ref{B2}}}{\leq} & \displaystyle\frac{5}{21}\cdot \frac{32}{43} + \frac{11}{43} = \frac{391}{903} < \frac{19}{43} \; ,
    \end{array}$$
    which is a contradiction. Therefore, $k=1$ or $k=2$, and so $G$ has order $10$ or $20$, respectively.

Clearly, the only non-cyclic group of order $10$ is $\mathcal{C}_5 \rtimes_{\iota} \mathcal{C}_2 \simeq P^*(5,2)$, which is an $M$-group. Assume that $|G|=20$. Hence, $G=P \rtimes H$, where $P\simeq\mathcal{C}_5$ and $H\simeq\mathcal{C}_4$.  Write~\mbox{$P=\langle x \rangle$} and \mbox{$H=\langle y \rangle$.}  By assumption, we have \mbox{$C_H(P)< H$,} so either \mbox{$C_H(P) = \{1\}$} or \mbox{$C_H(P) = \langle y^2 \rangle$.} If we suppose \mbox{$C_H(P) = \{1\}$,} then, by using the equality in Lemma~\ref{BoundPsiPNormalCyclicSylow}\ref{BoundPsiPNormalCyclicSylowFormulaSDP}, after dividing both sides by $\psi(\mathcal{C}_{|G|})$, we obtain
$$
\displaystyle \psi'(G)  = \frac{|P|}{\psi(\mathcal{C}_{|P|})}\cdot \psi'(H) + \left(1-\frac{|P|}{\psi(\mathcal{C}_{|P|})} \right)\cdot \frac{\psi(C_H(P))}{\psi(\mathcal{C}_{|H|})} = \frac{5}{21} + \left(1 - \frac{5}{21}\right) \cdot \frac{1}{11} = \frac{71}{231} < \frac{19}{43} \ ,$$
which is a contradiction. Therefore, $C_H(P) = \langle y^2 \rangle$. Then $x^{y^2} = x$, and it easily follows that $x^y = x^{-1}$. Thus, $G\simeq\mathcal{C}_5 \rtimes_{\iota} \mathcal{C}_4\simeq P^*(5,4)$, which is an $M$-group.

\bigskip

\noindent {\bf Case 2:} {\it $H\simeq \left(\mathcal{C}_{2^{k-1}}\!\times \mathcal{C}_2\right) \times \mathcal{C}_m$, where $k\geq 2$ and $(m,2)=1$.}

\smallskip

The maximal subgroups of $\mathcal{C}_{2^{k-1}}\!\times \mathcal{C}_2$ have order $2^{k-1}$, and since one of them is isomorphic to $\mathcal{C}_{2^{k-1}}$, it follows that $l(\mathcal{C}_{2^{k-1}}\!\times \mathcal{C}_2)=\psi(\mathcal{C}_{2^{k-1}})$. Thus, Lemma \ref{tech1}\ref{tech1.2} implies $l(H)=\psi(\mathcal{C}_{2^{k-1}} \times \mathcal{C}_m)$. Moreover, as observed in the proof of Proposition~\ref{2g-classification},
$$\psi'(H) = \frac{2^{2k}+5}{2^{2k+1}+1} \; .$$
Therefore,
    $$\begin{array}{ccl}
    \displaystyle \psi'(G) & \displaystyle \overset{\eqref{stareq}}{\leq} & \displaystyle\frac{|P|}{\psi(\mathcal{C}_{|P|})}\cdot\frac{2^{2k}+5}{2^{2k+1}+1} + \left(1-\frac{|P|}{\psi(\mathcal{C}_{|P|})} \right)\cdot\frac{2^{2k-1}+1}{2^{2k+1}+1}  \\[20pt]
    \displaystyle & \displaystyle = & \displaystyle\frac{|P|} {\psi(\mathcal{C}_{|P|})} \cdot \frac{2^{2k}-2^{2k-1}+4}{2^{2k+1}+1} + \frac{2^{2k-1}+1}{2^{2k+1}+1} \\[20pt]
    \displaystyle\ & \displaystyle \overset{k=2}{\leq} & \displaystyle \frac{|P|}{\psi(\mathcal{C}_{|P|})}\cdot \frac{4}{11} + \frac{3}{11}  \; \overset{\ref{BoundP/psi(C_P)}\text{\ref{B1}}}{\leq} \; \displaystyle\frac{3}{7}\cdot \frac{4}{11} + \frac{3}{11} = \frac{3}{7} < \frac{19}{43} \; .
    \end{array}$$

\bigskip

\noindent {\bf Case 3:} {\it $H\simeq M(2^k) \times \mathcal{C}_m$, where $k\geq 4$ and $(m,2) =1$.}

\smallskip

 In this case, the analysis is exactly the same as in the previous one, in view of the fact that \hbox{$l(M(2^k))=\psi(\mathcal{C}_{2^{k-1}})$,} so again $l(H)=\psi(\mathcal{C}_{2^{k-1}} \times \mathcal{C}_m)$ (by Lemma~\ref{tech1}\ref{tech1.2}), and that $\psi'(M(2^k) \times \mathcal{C}_m)=\psi'(\left(\mathcal{C}_{2^{k-1}} \!\times \mathcal{C}_2\right) \times \mathcal{C}_m)$ (by direct computation, as in the proof of Proposition~\ref{2g-classification}). Note that the only difference lies in the domain of $k$, which is now $k\geq 4$, whereas previously it was $k\geq 2$: this does not affect the validity of the previously established upper bound, since it was obtained by evaluating a function that is decreasing in the variable $k$ at $k=2$.

\bigskip

\noindent {\bf Case 4:} {\it $H\simeq Q_8 \times \mathcal{C}_m$, where $(m,2) =1$.}

\smallskip

The maximal subgroups of $Q_8$ are all isomorphic to $\mathcal{C}_4$, so $l(Q_8)=\psi(\mathcal C_4)$ and consequently, by~Lem\-ma~\ref{tech1}\ref{tech1.2}, it follows $l(H)=\psi(\mathcal{C}_4\times\mathcal{C}_m)$. Moreover, $\psi'(H)=\frac{27}{43}$. Therefore,
$$\begin{array}{ccl}
\displaystyle \psi'(G) & \displaystyle \overset{\eqref{stareq}}{\leq} & \displaystyle\frac{|P|}{\psi(\mathcal{C}_{|P|})}\cdot\frac{27}{43} + \left(1-\frac{|P|}{\psi(\mathcal{C}_{|P|})} \right)\cdot\frac{11}{43} = \displaystyle\frac{|P|}{\psi(\mathcal{C}_{|P|})} \cdot \frac{16}{43} + \frac{11}{43} \\[20pt]
\displaystyle\ & \displaystyle \overset{\text{\ref{BoundP/psi(C_P)}\ref{B1}}}{\leq} & \displaystyle\frac{3}{7}\cdot \frac{16}{43} + \frac{11}{43} = \frac{125}{301} < \frac{19}{43} \; .
\end{array}$$

\bigskip

\noindent {\bf Case 5:} {\it $H \simeq \left(\mathcal{C}_3 \rtimes_{\iota} \mathcal{C}_{2^k}\right) \times \mathcal{C}_m$, where $k\geq 1$ and $(m,6)=1$.}

\smallskip

It can be readily shown that the maximal subgroups of~$\mathcal{C}_3 \rtimes_{\iota} \mathcal{C}_{2^k}$ are isomorphic either to~$\mathcal{C}_{2^k}$ or to $\mathcal{C}_{3\cdot 2^{k-1}}$. Hence, since $\psi(\mathcal{C}_{2^k}) < \psi(\mathcal{C}_{3\cdot 2^{k-1}})$ for every $k\geq 1$ (by direct computation), we have $l(\mathcal{C}_3 \rtimes_{\iota} \mathcal{C}_{2^k}) = \psi(\mathcal{C}_{3\cdot 2^{k-1}})$, and so Lem\-ma~\ref{tech1}\ref{tech1.2} yields $l(H)=\psi(\mathcal{C}_{3\cdot 2^{k-1}}\times \mathcal{C}_m)$. Moreover, using Proposition~\ref{PsiP*}, we obtain
$$
\displaystyle\psi'\big(\mkern-1mu \left(\mathcal{C}_3 \rtimes_{\iota} \mathcal{C}_{2^k}\right) \times \mathcal{C}_m {\mkern1.5mu} \big) = \frac{2^{2k+3}+7}{7(2^{2k+1} + 1)} \; .$$
Therefore, also taking into account that in this case $p\geq 5$, we have
$$\begin{array}{ccl}
    \displaystyle \psi'(G) & \displaystyle \overset{\eqref{stareq}}{\leq} & \displaystyle\frac{|P|}{\psi(\mathcal{C}_{|P|})}\cdot\frac{2^{2k+3}+7}{7(2^{2k+1} + 1)} + \left(1-\frac{|P|}{\psi(\mathcal{C}_{|P|})} \right)\cdot\frac{2^{2k-1}+1}{2^{2k+1}+1} \;\\[20pt] 
    \displaystyle & \displaystyle = & \displaystyle\frac{|P|} {\psi(\mathcal{C}_{|P|})} \cdot \frac{2^{2k+3}-7\cdot 2^{2k-1}}{7(2^{2k+1}+1)} + \frac{2^{2k-1}+1}{2^{2k+1}+1} \; \\[20pt]
    \displaystyle & \displaystyle \leq & \displaystyle\frac{|P|} {\psi(\mathcal{C}_{|P|})} \, \cdot \, \lim_{k\to +\infty}\Bigg[ \frac{2^{2k+3}-7\cdot 2^{2k-1}}{7(2^{2k+1}+1)} \Bigg] + \, \frac{2^{2k-1}+1}{2^{2k+1}+1} \; \\[20pt]
    \displaystyle & \displaystyle = & \displaystyle\frac{|P|} {\psi(\mathcal{C}_{|P|})} \cdot \frac{9}{28} + \frac{2^{2k-1}+1}{2^{2k+1}+1} \; \overset{k=1}{\leq} \; \displaystyle\frac{|P|} {\psi(\mathcal{C}_{|P|})} \cdot \frac{9}{28} + \frac{1}{3} \; \\[20pt]
    \displaystyle\ & \displaystyle \overset{\text{\ref{BoundP/psi(C_P)}\ref{B2}}}{\leq} & \displaystyle\frac{5}{21}\cdot \frac{9}{28} + \frac{1}{3} = \frac{241}{588} < \frac{19}{43} \; . \end{array}$$

\bigskip

\noindent {\bf Case 6:} {\it $H \simeq \left(\mathcal{C}_5 \rtimes_{\iota} \mathcal{C}_2\right) \times \mathcal{C}_m$, where $(m,10)=1$.}

\smallskip

The maximal subgroups of $\mathcal{C}_5 \rtimes_{\iota} \mathcal{C}_2$ are isomorphic either to $\mathcal{C}_2$ or to $\mathcal{C}_5$, and so \mbox{$l(\mathcal{C}_5 \rtimes_{\iota} \mathcal{C}_2)=\psi(\mathcal{C}_5)$.} Thus, by Lemma \ref{tech1}\ref{tech1.2}, $l(H)=\psi(\mathcal{C}_5\times \mathcal{C}_m)$. Moreover, using Proposition \ref{PsiP*}, we find that $\psi'(H) = \frac{31}{63}$. Therefore,
$$\begin{array}{ccl}
\displaystyle \psi'(G) & \displaystyle \overset{\eqref{stareq}}{\leq} & \displaystyle\frac{|P|}{\psi(\mathcal{C}_{|P|})}\cdot\frac{31}{63} + \left(1-\frac{|P|}{\psi(\mathcal{C}_{|P|})} \right)\cdot\frac{1}{3} = \displaystyle\frac{|P|}{\psi(\mathcal{C}_{|P|})} \cdot \frac{10}{63} + \frac{1}{3}\\[20pt]
\displaystyle\ & \displaystyle \overset{\text{\ref{BoundP/psi(C_P)}\ref{B1}}}{\leq} & \displaystyle\frac{3}{7}\cdot \frac{10}{63} + \frac{1}{3} = \frac{59}{147} < \frac{19}{43} \; .
\end{array}$$

\bigskip

\noindent {\bf Case 7:} {\it $H \simeq \left(\mathcal{C}_5 \rtimes_{\iota} \mathcal{C}_4\right) \times \mathcal{C}_m$, where $(m,10)=1$.}

\smallskip

The maximal subgroups of $\mathcal{C}_5 \rtimes_{\iota} \mathcal{C}_4$ are isomorphic either to $\mathcal{C}_4$ or to $\mathcal{C}_{10}$, and so \hbox{$l(\mathcal{C}_5\!\rtimes_{\iota} \!\mathcal{C}_4)\!=\!\psi(\mathcal{C}_{10})$}. Thus, by Lemma \ref{tech1}\ref{tech1.2}, $l(H) = \psi(\mathcal{C}_{10} \times \mathcal{C}_m)$. Moreover, by~Proposition~\ref{PsiP*}, $\psi'(H) = \frac{103}{231}$. Therefore,
$$\begin{array}{ccl}
\displaystyle \psi'(G) & \displaystyle\overset{\eqref{stareq}}{\leq} & \displaystyle\frac{|P|}{\psi(\mathcal{C}_{|P|})}\cdot\frac{103}{231} + \left(1-\frac{|P|}{\psi(\mathcal{C}_{|P|})} \right)\cdot\frac{3}{11} = \displaystyle\frac{|P|}{\psi(\mathcal{C}_{|P|})} \cdot \frac{40}{231} + \frac{3}{11} \\[20pt]
\displaystyle\ & \displaystyle \overset{\text{\ref{BoundP/psi(C_P)}\ref{B1}}}{\leq} & \displaystyle\frac{3}{7}\cdot \frac{40}{231} + \frac{3}{11} = \frac{17}{49} < \frac{19}{43} \; .
\end{array}$$

The proof that $G$ is an $M$-group is now complete. Finally, by applying Proposition~\ref{Mg1943} once more, we obtain the stated classification.  \hfill $\square$

\begin{remark}\label{changeinduction}
    The inductive argument in the proof of Theorem A proceeds by showing that any group with the given bound on $\psi'$ is an M-group. In this context, to apply the inductive hypothesis, it is first necessary to identify all groups satisfying the criterion being established. This information is provided by Proposition~\ref{Mg1943}, which relies on knowledge of the structure of finite M-groups. However, once the list of groups satisfying the criterion is known, knowledge of M-groups is not really needed to prove Theorem A. Indeed, it suffices to employ an inductive argument showing that any group with the given bound on $\psi'$ belongs to that list, so that the inductive hypothesis already incorporates the information provided by Proposition~\ref{Mg1943}. Proceeding in this alternative way for~The\-o\-rem~A, all possible groups appearing in the conclusion arise when $[P,H]=\{1\}$ or in the non-contradictory instances of Case 1. In light of this, Proposition~\ref{Mg1943} should just be regarded as an elegant method --- one that exploits the known structure of finite M-groups --- to identify the groups satisfying the criterion, thereby avoiding the use of computational algebra software.
    
This observation will be useful in the final considerations of the paper, where we outline the path leading to the classification of all groups whose $\psi'$-value exceeds $\frac{31}{77}$.
\end{remark}

\ \\
\noindent {\bf Proof of Theorem B.}\\
Clearly $\psi'\big(D_8\times\mathcal{C}_m\big) = \frac{19}{43}$, and, by Proposition \ref{PsiP*}, $\psi'\big((\mathcal{C}_7 \rtimes_\iota \mathcal{C}_2) \times \mathcal{C}_m\big) = \frac{19}{43}$.
We proceed by induction on the order of $G$ to show that these are the only groups with $\psi'$-value equal to $\frac{19}{43}$.
Clearly,~$G$ is not cyclic, and so $2$ divides $|G|$ by Lemma \ref{smallestprime1943-3177}\ref{smallestprime1943}.
If $G$ is a $2$-group, then by Lemma~\ref{2g-classification}\ref{2g=1943} we have $G \simeq D_8$, and the claim follows. We may thus assume that~$G$ is not a $2$-group. Therefore, Corollary~\ref{1943SylNormCycl} yields that
$$G=P\rtimes H \; ,$$
where $P$ is a cyclic Sylow $r$-subgroup of $G$ with $r=2$ or $r=p\coloneqq\max\pi(G)$, and $H$ is a suitable subgroup of $G$. If $G=P\times H$, then $\psi'(H) = \frac{19}{43}$, so, by the induction hypothesis, $H$ is one of the groups listed in the statement and the claim follows immediately. Therefore, we can assume that $[P,H] \neq \{1\}$. Consequently, $r=p$ (otherwise $r=2$, and so~\mbox{$[P,H]=\{1\}$).}

By Lemma \ref{BoundPsiQuotient},
$$\frac{19}{43} = \psi'(G) \leq \psi'(P)\psi'(G/P)=\psi'(H) \, ,$$
and since $P$ is not central in $G$, the equality does not hold. Consequently, $\psi'(H)>\frac{19}{43}$ and so $H$ is one of the groups classified by Theorem A.

From this point on, the proof proceeds as in Theorem A, with some modifications that we will now point out and analyze. To ensure clarity and consistency, both proofs are presented using the same notation. 

If $H$ is non-cyclic, the proof remains unchanged: we always obtain the contradiction $\psi'(G)<\frac{19}{43}$. Therefore, assume that $H$ is cyclic. Set $|P|=p^n$ and let $k$ be the positive integer such that $2^k$ is the highest power of $2$ dividing $|H|$. First, if we suppose that $n\geq 2$, then we obtain a contradiction just as before, so $n=1$. Next, the assumption $p\geq 7$ leads to~\hbox{$\psi'(G)\leq\frac{19}{43}$,} which is now no longer a contradiction. Instead, the assumption $p\geq 11$ leads to the contradiction
$$\begin{array}{c}
    \displaystyle \psi'(G) \displaystyle \leq  \displaystyle\frac{|P|} {\psi(\mathcal{C}_{|P|})}\cdot\frac{2}{3} + \frac{1}{3}\overset{\text{\ref{BoundP/psi(C_P)}\ref{Bplus}}}{\leq} \displaystyle\frac{11}{111}\cdot \frac{2}{3} + \frac{1}{3} = \frac{133}{333} < \frac{19}{43} \; .   
    \end{array}$$

\noindent Consequently, $p$ can be equal to $3$, $5$ or $7$.

Assuming $|P|=3$ again leads to the conclusion that $G$ is isomorphic to $P^*(3,2^k)$, thereby yielding a contradiction according to Corollary~\ref{P*(2)1943}.

Proceeding under the assumption that $|P| = 5$, it follows once again that $G$ must be isomorphic to $P^*(5,2)$ or $P^*(5,4)$. However, by Corollary \ref{P*(2)1943}, these groups have $\psi'$-value exceeding $\frac{19}{43}$, leading to a contradiction.

Finally, we have to deal with the case $|P|=7$ (note that the case $p=7$ did not require a separate treatment in the proof of Theorem A). Since $p=\max\pi(G)$, we have $\pi(G) \subseteq \{2,3,5\}$, so $H=H_2\times H_3 \times H_5$, where each $H_i$ denotes the Sylow $i$-subgroup of $H$. Given that \hbox{$\operatorname{Aut}(P) \simeq \mathcal{C}_6$,} it follows that $[P,H_5]=\{1\}$, so $H_5\leq Z(G)$; in particular, \hbox{$G=H_5\times \big(P \rtimes (H_2\times H_3) \big)$.} If~\hbox{$H_5 \neq\{1\}$,} then $\psi'\big(P \rtimes (H_2\times H_3)\big)=\frac{19}{43}$, and the claim follows by the induction hypothesis. Assume therefore that $H_5=\{1\}$, so~\hbox{$G = P \rtimes (H_2\times H_3)$.}
Since $2$ divides $|G|$, $H_2\neq \{1\}$. If $H_3 = \{1\}$, then $G\simeq P^*(7,2^k)$. Therefore, by~Pro\-position~\ref{PsiP*},
    $$\psi'(G) = \frac{2^{2k+5}+43}{43\cdot2^{2k+1}+43} \; , $$
which is readily seen to equal $\frac{19}{43}$ if and only if $k=1$. Hence, $G\simeq \mathcal{C}_7 \rtimes_{\iota} C_2$ and the claim follows. Thus, we may further assume that $H_3\neq\{1\}$.

Now, if $P$ centralizes one of the $H_i$'s, then the argument can be concluded inductively, as in the case $H_5 \neq \{1\}$. Therefore, let us assume that $P$ centralizes neither~$H_2$ nor~$H_3$. Set $|H|=2^k\cdot 3^t$. Since $H/C_H(P)$ embeds in $\operatorname{Aut}(P) \simeq \mathcal{C}_6$ and \hbox{$H_2,H_3 \not\leq C_H(P)$,} it follows that $|C_H(P)| = 2^{k-1}\cdot 3^{t-1}$. Therefore,
    $$\begin{array}{ccl}
    \displaystyle \psi'(G) & \displaystyle\overset{\text{\ref{BoundPsiPNormalCyclicSylow}\ref{BoundPsiPNormalCyclicSylowFormulaSDP}}}{\leq} & \displaystyle\frac{|P|}{\psi(\mathcal{C}_{|P|})} \cdot 1 + \left(1-\frac{|P|}{\psi(\mathcal{C}_{|P|})} \right)\cdot\frac{2^{2k-1}+1}{2^{2k+1}+1} \cdot \frac{3^{2t-1}+1}{3^{2t+1}+1}  \\[20pt]
    \displaystyle & \displaystyle \overset{k=1}{\leq}  \displaystyle  & \displaystyle\frac{|P|}{\psi(\mathcal{C}_{|P|})} + \left(1-\frac{|P|}{\psi(\mathcal{C}_{|P|})} \right)\cdot\frac{1}{3} \cdot \frac{3^{2t-1}+1}{3^{2t+1}+1}\\[20pt]
    \displaystyle & \displaystyle \overset{t=1}{\leq} & \displaystyle \frac{|P|}{\psi(\mathcal{C}_{|P|})} + \left(1-\frac{|P|}{\psi(\mathcal{C}_{|P|})} \right)\cdot\frac{1}{3} \cdot \frac{1}{7} \, = \,
    \displaystyle\frac{|P|} {\psi(\mathcal{C}_{|P|})}\cdot\frac{20}{21} + \frac{1}{21} \\[20pt] 
    \displaystyle & \displaystyle \overset{\text{\ref{BoundP/psi(C_P)}\ref{B3}}}{\leq} & \displaystyle\frac{7}{43}\cdot \frac{20}{21} + \frac{1}{21} = \frac{61}{301} < \frac{19}{43} \; , \end{array}$$
which is a contradiction. This concludes the proof. \hfill $\square$

\bigskip
\bigskip

Finally, we wish to briefly explain to the reader how to deal with groups whose \hbox{$\psi'$-value} lies in the interval $\big(\frac{31}{77},1\big]$ (these groups are listed in the introduction). The proof scheme essentially follows that of Theorem~A: we express $G$ as a semidirect product $P\rtimes H$ (here $|P|=p^n$, where $p\geq3$ is the largest prime dividing $|G|$, while $|H|$ is even and $2^k$ is the largest power of $2$ dividing $|H|$), and then apply induction to obtain that $H$ is one of the required groups. The proof is subsequently divided into cases, accordingly to the form of~$H$ (note that the thesis is obvious if $[P,H]=\{1\}$, so one assumes that $[P,H]\neq\{1\}$), and the arguments still rely on inequality~\eqref{stareq}.    

If $H$ is cyclic, then the same computation as in the corresponding case of the proof of Theorem A yields $|P|\in\{3,5,7,9\}$. Then we can further assume that $|H|=2^k3^m$, with $m=0$ when $p\in\{3,5, 9\}$.  If $|P|=3$, then $G\simeq \mathcal C_3\rtimes_\iota \mathcal C_{2^k}$. If $|P|=9$, then the same computation with $k=2$ gives a contradiction, so $G\simeq D_{18}$. If $|P|=5$, a direct computation shows that the action of $H$ on $P$ is by inversion, so $G\simeq \mathcal C_5\rtimes_\iota\mathcal C_{2^k}$. If $|P|=7$, then we deduce as above that $k=1$, and thereafter we directly compute the $\psi'$-value of~\hbox{$G\simeq\mathcal C_7\rtimes\big(\mathcal C_2\times\mathcal C_{3^m}\big)$,} obtaining $m=0$ and $G\simeq D_{14}$.

If~$H\simeq\mathcal C_3\times\mathcal C_3\times \mathcal C_m$ with $(m,3)=1$, then we have \mbox{$\psi'(H) =$ \kern0.12em\scalebox{0.8}{$\dfrac{25}{61}$}} and, as a consequence of Lemma~\ref{techCycl},~\scalebox{0.8}{$\dfrac{l(H)}{\psi(\mathcal{C}_{|H|})} \leq \dfrac{25}{183}$}, which combined with inequality~\eqref{stareq} leads to a contradiction. If $H\simeq \mathcal C_{2^{k-1}}\!\times \mathcal C_2\times \mathcal C_m$ with $(m,2)=1$ and $k\geq2$, we first obtain $|P|=3$ and then $k=2$, which gives $G\simeq D_{12}\times \mathcal C_m$. If $H\simeq M(2^k)\times\mathcal C_m$, where $k\geq4$ and $(m,2)=1$, then we obtain a contradiction. If~\mbox{$H\simeq Q_{16}\times\mathcal C_m$} with $(m,2)=1$, then we have \mbox{$\psi'(H) =$ \kern0.12em\scalebox{0.8}{$\dfrac{75}{171}$}} (by direct computation) and $l(H)=\psi(\mathcal{C}_8)$ (by Lemma~\ref{tech1}\ref{tech1.2}), so inequality~\eqref{stareq} yields a contradiction. If~\mbox{$H\simeq Q_8\times\mathcal C_m$} with $(m,2)=1$, then we obtain $p=3$ and subsequently $n=1$, so $G\simeq \big(\mathcal C_3\rtimes Q_8\big)\times\mathcal C_m$.

Suppose $H$ is non-abelian and not a primary group. In this case, $H=(H'\rtimes D)\times C$, where~$D$ is either cyclic of order a power of $2$ or $Q_8$, $H'$ is cyclic with $|H'|\in \{3,5,7,9\}$, and $C$ is cyclic with $(|C|,|H'DP|)=1$. Since $H/C_H(P)$ is abelian, we have that $C_H(P)$ contains $H'$. Then $H'$ is normal in $G$ and~$G/H'$ must be one of the required groups. By symmetry, we may thus assume that \hbox{$C_H(P)\neq H$} and $|P|\in \{3,5,7,9\}$ with $(|H'|,|P|)=1$. This in particular implies that $D$ is cyclic, so~\hbox{$G\simeq C\times \big((P\times H')\rtimes_\iota \mathcal C_{2^k}\big)$}. If~\hbox{$k\geq3$,} then we may further assume $|P|=5$ and $|H'|=3$, and we can easily rule out this case by direct computation of the $\psi'$. If $k\leq2$, then there is only a finite number of possibilities, which can either be ruled out by GAP or by direct computation.

\bigskip

\paragraph{Acknowledgments}

The authors are members of the non-profit association ‘‘AGTA --- Advances in Group Theory and Applications’’ (www.advgrouptheory.com), and are supported by GNSAGA (INdAM). Funded by the European Union - Next Generation EU, Missione 4 Componente 1 CUP B53D23009410006, PRIN 2022- 2022PSTWLB - Group Theory and Applications.

\bigskip
\bigskip

\begin{flushleft}
\rule{8cm}{0.4pt}\\
\end{flushleft}

\bigskip
\bigskip

{
\sloppy
\noindent
Luigi Iorio

\noindent 
Dipartimento di Matematica e Applicazioni ``Renato Caccioppoli''

\noindent
Università degli Studi di Napoli Federico II

\noindent
Complesso Universitario Monte S. Angelo

\noindent
Via Cintia, Napoli (Italy)

\noindent
e-mail: luigi.iorio2@unina.it 

}

\bigskip
\bigskip

{
\sloppy
\noindent
Marco Trombetti

\noindent 
Dipartimento di Matematica e Applicazioni ``Renato Caccioppoli''

\noindent
Università degli Studi di Napoli Federico II

\noindent
Complesso Universitario Monte S. Angelo

\noindent
Via Cintia, Napoli (Italy)

\noindent
e-mail: marco.trombetti@unina.it 

}


\begin{thebibliography}{99}

\bibitem{FIRST} H. Amiri, S.M. Jafarian Amiri, and I.M. Isaacs: ‘‘Sums of element orders in finite groups’’, {\it Comm. Algebra} 37 (2009), 2978--2980.

\bibitem{SolvableCriterion} M. Baniasad Azad and B. Khosravi:
‘‘A criterion for solvability of a finite group by the sum of element orders’’, {\it J. Algebra} 516 (2018), 115--124.

\bibitem{SSCriterion} M. Baniasad Azad and B. Khosravi: ‘‘On two conjectures about the sum of element orders’’, {\it Canad. Math. Bull.} 65 (2022), 30--38.

\bibitem{uguaglianza} A. Bahri, B. Khosravi and Z. Akhlaghi:
‘‘A result on the sum of element orders of a finite group’’,
{\it Arch. Math. \textnormal(Basel\textnormal)} 114 (2020), no. 1, 3–12.


\bibitem{GAP} The GAP Group: ‘‘GAP – Groups, Algorithms, and Programming’’, v4.14.0 (2024).

\bibitem{Gorenstein} D. Gorenstein: ‘‘Finite Groups’’, 2nd ed., \textit{Chelsea}, New York (1980).

\bibitem{NonCyclicCriterionI} M. Herzog, P. Longobardi, and M. Maj: ‘‘An exact upper bound for sums of element orders in non-cyclic finite groups’’, {\it J. Pure Appl. Algebra} 222 (2018), 1628--1642.

\bibitem{GeneralizedBoundPsiCyclic} M. Herzog, P. Longobardi, and M. Maj: ‘‘Two new criteria for solvability of finite groups’’, {\it J. Algebra} 511 (2018), 215--226.

\bibitem{NonCyclicCriterionII} M. Herzog, P. Longobardi, and M. Maj: ‘‘The second maximal groups with respect to the sum of element orders’’, {\it J. Pure Appl. Algebra} 225 (2021), 106531, 11pp.

\bibitem{Lucchini} A. Lucchini: ‘‘On the order of transitive permutation groups with cyclic point-stabilizer’’, {\it Rend. Lincei, Mat. Appl.} 9 (1998), 241--243.
 

\bibitem{Schmidt} R. Schmidt: ‘‘Subgroup Lattices of Groups’’, {\it de Gruyter}, Berlin (1994).

\bibitem{SuzukiBooks} M. Suzuki: ‘‘Group Theory I, II’’, {\it Springer}, Berlin (1986).

\bibitem{NilpotentCriterion} M. T\u{a}rn\u{a}uceanu: ‘‘A criterion for nilpotency of a finite group by the sum of element orders’’, {\it Comm. Algebra} 49 (2021), no. 4, 1571--1577.



\end{thebibliography}
\end{document}